\documentclass[12pt]{article}
\usepackage{latexsym}
\usepackage[dvipdfmx]{graphicx}
\newsavebox{\toy}
\savebox{\toy}{\framebox[0.65em]{\rule{0cm}{1ex}}}
\newcommand{\QED}{\usebox{\toy}\end{demo}}
\newenvironment{property}%
{\begin{list}{}{\setlength{\rightmargin}{0pt}%
\setlength{\itemsep}{0pt}}}{\end{list}}
\newlength{\templength}
\newcommand{\bp}{\setlength{\templength}{\labelwidth}%
\setlength{\labelwidth}{2em}\begin{property}}
\newcommand{\ep}{\end{property}\setlength{\labelwidth}{\templength}}
\newtheorem{theorem}{Theorem}[subsection]
\newtheorem{lemma}[theorem]{Lemma}
\newtheorem{proposition}[theorem]{Proposition}
\newtheorem{corollary}[theorem]{Corollary}
\newtheorem{definition}[theorem]{Definition}
\newtheorem{example}[theorem]{Example}
\newtheorem{remark}[theorem]{Remark}
\newtheorem{exercise}{Exercise}[subsection]
\newtheorem{assumption}{Assumption}
\newcommand{\Thm}[1]{Theorem \ref{Thm.#1}}
\newcommand{\Lem}[1]{Lemma \ref{Lem.#1}}

\newcommand{\Prop}[1]{Proposition \ref{Prop.#1}}
\newcommand{\Cor}[1]{Corollary \ref{Cor.#1}}

\newcommand{\Theorem}[1]{\begin{theorem}\label{Thm.#1}}
\newcommand{\stoptheorem}{\end{theorem}}
\newcommand{\Lemma}[1]{\begin{lemma}\label{Lem.#1}}
\newcommand{\stoplemma}{\end{lemma}}
\newcommand{\Proposition}[1]{\begin{proposition}\label{Prop.#1}}
\newcommand{\stopproposition}{\end{proposition}}
\newcommand{\Corollary}[1]{\begin{corollary}\label{Cor.#1}}
\newcommand{\stopcorollary}{\end{corollary}}
\newcommand{\Assumption}[1]{\begin{assumption}\label{Ass.#1}\rm}
\newcommand{\Definition}[1]{\begin{definition}\label{Def.#1}\rm}
\newcommand{\stopdefinition}{\end{definition}}
\newcommand{\Remark}[1]{\begin{remark}\label{Rem.#1}\rm }
\newcommand{\Exercise}[1]{\begin{exercise}\label{Exe.#1}\rm }
\newcommand{\Example}[1]{\begin{example}\label{Exa.#1}\rm }
\newcommand{\stopexample}{\end{example}}
\newcommand{\bd}{\begin{displaymath}}
\newcommand{\ed}{\end{displaymath}}
\newcommand{\bdn}{\begin{equation}}
\newcommand{\bdnl}{\begin{equation}\label}
\newcommand{\edn}{\end{equation}}
\newcommand{\barray}{\begin{array}}
\newcommand{\earray}{\end{array}}
\newcommand{\bds}{\begin{description}}
\newcommand{\eds}{\end{description}}
\newcommand{\bitemize}{\begin{itemize}}
\newcommand{\eitemize}{\end{itemize}}
\newcommand{\benumerate}{\begin{enumerate}}
\newcommand{\eenumerate}{\end{enumerate}}
\newcommand{\btabbing}{\begin{tabbing}}
\newcommand{\etabbing}{\end{tabbing}}
\newcommand{\bcenter}{\begin{center}}
\newcommand{\ecenter}{\end{center}}
\newcommand{\bflushright}{\begin{flushright}}
\newcommand{\bflushleft}{\begin{flushleft}}
\newcommand{\eflushright}{\end{flushright}}
\newcommand{\eflushleft}{\end{flushleft}}
\newcommand{\bdnn }{\begin{eqnarray*}}
\newcommand{\ednn }{\end{eqnarray*}}
\newcommand{\bdmn}{\begin{eqnarray}}
\newcommand{\edmn}{\end{eqnarray}}

\newcommand{\nn}{\nonumber}
\newcommand{\SSC}[1]{\section{#1}\setcounter{equation}{0}}

\newcounter{biblio}
\newenvironment{references}%
{\begin{list}{[\arabic{biblio}]}{\usecounter{biblio}%
\setlength{\leftmargin}{2.5em}\setlength{\rightmargin}{0pt}%
\setlength{\labelwidth}{2em}\setlength{\itemsep}{0pt}}}{\end{list}}
\newcommand{\References}%
{\vspace{2.8ex plus .3ex minus .3ex}%
\begin{center}{\bf References}\end{center}\begin{references}}

\usepackage{amssymb}
\usepackage{mathrsfs}
\usepackage{amsmath}
\usepackage{array}
\usepackage{makeidx}
\makeindex

\usepackage{ascmac}
\usepackage{wrapfig}
\newcommand{\N}{{\mathbb{N}}}
\newcommand{\Z}{{\mathbb{Z}}}
\newcommand{\zd}{\Z^d}
\newcommand{\Q}{{\mathbb{Q}}}
\newcommand{\R}{{\mathbb{R}}}
\newcommand{\rd}{\R^d}
\newcommand{\C}{{\mathbb{C}}}

\newcommand{\ra }{\rightarrow }
\newcommand{\lra }{\longrightarrow }

\newcommand{\Ra}{\Rightarrow }
\newcommand{\La}{\Leftarrow }
\newcommand{\Lra}{\Longrightarrow }
\newcommand{\Llra}{\Longleftrightarrow }

\newcommand{\ov}{\overline}
\newcommand{\tl}{\widetilde}

\newcommand{\vvs}{\vspace{2ex}}
\newcommand{\vs}{\vspace{1ex}}

\newcommand{\Ran}{{\rm Ran }}
\newcommand{\Ker}{{\rm Ker }}
\newcommand{\lan}{\langle \:}
\newcommand{\ran}{\: \rangle}
\newcommand{\lef}{\left}

\newcommand{\ri}{\right}
\newcommand{\st}{\stackrel}

\newcommand{\8}{\infty}

\newcommand{\6}{\partial}

\newcommand{\sub}{\subset}

\newcommand{\bsh}{\backslash}

\renewcommand{\Re}{\mathop{\rm Re}\nolimits}

\newcommand{\inflim}{\mathop{\underline{\lim}}}

\newcommand{\epty}{\emptyset}
\renewcommand{\a}{\alpha}
\renewcommand{\b}{\beta}
\newcommand{\gm}{\gamma}

\newcommand{\del}{\delta}
\newcommand{\D}{\Delta}

\newcommand{\e}{\varepsilon}

\newcommand{\tht}{\theta}

\newcommand{\lm}{\lambda}
\newcommand{\Lm}{\Lambda}

\newcommand{\m}{\mu}

\newcommand{\rh}{\rho}

\newcommand{\s}{\sigma}

\newcommand{\vp}{\varphi}

\newcommand{\w}{\omega}

\newcommand{\W}{\Omega}

\newcommand{\cA }{{\cal A}}
\newcommand{\cB }{{\cal B}}

\newcommand{\cF }{{\cal F}}

\newcommand{\bfi}{{\bf i}}

\newcommand{\deq}{\st{\rm def}{=}}
\newcommand{\supp}{{\rm supp\:}}
\newcommand{\vertiii}[1]{{\left\vert\kern-0.25ex\left\vert\kern-0.25ex\left\vert #1 
    \right\vert\kern-0.25ex\right\vert\kern-0.25ex\right\vert}}
\makeatletter
\def\section{\@startsection{section}{1}{\z@}{-3.5ex plus -1ex minus 
 -.2ex}{2.3ex plus .2ex}{\bf}}
\def\subsection{\@startsection{subsection}{2}{\z@}{-3.25ex plus -1ex minus 
 -.2ex}{1.5ex plus .2ex}{\bf}}
\makeatother

\begin{document}

\bcenter
\large{\bf The transition operator of a
  random walk perturbated by sparse potentials}\footnote{private notes \today. }

\normalsize 

\vvs

\vvs Takuya Mine\footnote{Faculty of Arts and Sciences,
  Kyoto Institute of Technology,
  Matsugasaki, Sakyo-ku, Kyoto
606-8585, Japan} 
and Nobuo Yoshida\footnote{Graduate School of Mathematics,
  Nagoya University,
Furocho, Chikusaku, Nagoya 464-8602, Japan}

\ecenter
\begin{abstract}
  We consider an operator $P_V=(1+V)P$ on $\ell^2(\Z^d)$, where $P$
  is the transition operator 
  of a symmetric irreducible random walk,
  and $V$ is a ``sparse'' potential. We first characterize 
  the essential spectra of this operator. Secondly, we
  prove that all the eigenfunctions which correspond to discrete spectra
  decay exponentially fast. Thirdly, we give a sufficient
  condition for this operator to have an absolute
  spectral gap at the right edge of the spectra.
  Finally, as an application of the absolute spectral gap
  and the exonential decay of the eigenfunctions, 
  we prove a limit theorem for the random walk under the
  Gibbs measure associated to the potential $V$. 
\end{abstract}
\SSC{Introduction}
In this article, we investigate the spectral properties of
an operator 
  $P_V=(1+V)P$ on $\ell^2(\Z^d)$, where $P$
  is the transition operator 
  of a symmetric irreducible random walk,
  and $V$ is a nonnegative bounded function.
  Here, $V$ is supposed to be a so called ``sparse'' potential,
  of which a typical example is that with the property
  \bdnl{sparse*0}
  \min \{|x-y|\; ;\; x,y \in {\rm supp} V, \; x \neq y,\; |x| \ge r,\; |y| \ge r\}
  \st{r \ra \8}{\lra}\8,
  \edn 
    Before introducing the contents of this article, we start by explaining the probabilistic
  backgrounds which brought us to the study of the spectra of $P_V$.

  Let $(S_n)_{n \in \N}$ be a random walk on a probability space
  $(\W, \cF, P)$ which is associated with the transition operator $P$.
  We suppose that $S_0=0$. 
  Then, the semigroup $P_V^n$ ($n \in \N$)
  is expressed by the discrete version of the Feynmann-Kac formula.
\bdnl{FKPV}
P_V^nf(x)=E\lef[ f(x+S_n)\prod^{n-1}_{j=0}(1+V(x+S_j))\ri],
\; \; f \in \ell^2(\zd),\; \; x \in \Z^d. 
\edn
Similarly, if we restrict $P_V$ to
$\ell^2(\N^d)$ by imposing Dirichlet boundary condition,
and denote the restriction by $P^+_V \in \cB (\ell^2(\N^d))$,
then, we obtain
\bdmn
(P_V^+)^nf(x)
&=& E\lef[ f(x+S_n){\bf 1}\{ S_j \in \N^d,\; j=1,\ldots,n\}
  \prod^{n-1}_{j=0}(1+V(x+S_j))\ri],\nn \\
& &f \in \ell^2(\N^d),\; \; x \in \N^d.
\label{FKPV+}
\edmn
Let $\6 \N^d=\bigcup_{\a=1}^d\{x \in \N^d\; ; \; x_\a=0\}$, 
and $\b>0$ be a positive parameter. 
In connection with statistical physics, the following
sequence $\m_N$, $N \ge 1$ of measures on $(\W,\cF)$ are studied.
\bdnl{Pn+}
\m_N(d\w)={1 \over Z_N}E\lef[{\bf 1}\{ S_j \in \N^d,\; j=1,\ldots,N\}
  \exp \lef( \b \sum_{j=1}^{N-1}{\bf 1}_{\6 \N^d}(S_j)\ri) : d \w\ri], 
\edn 
where $Z_N$ is the normalizing constant. The above measure is related with the
operator $P_V^+$ with $V=\exp (\b{\bf 1}_{\6 \N^d})-1$ via the formula (\ref{FKPV+}). 
Under the measure $\m_N$, two competing effects coexist: The paths of the
random walk is attracted to $\6 \N^d$ by the potential to lower the eneregy:
$$
-\b \sum_{j=1}^{N-1}{\bf 1}_{\6 \N^d}(S_j). 
$$
On the other hand, paths of the
random walk are pushed away from $\6 \N^d$ by the entropic repulsion.
These competing effects cause the
phase transition, called wetting transition, which is expressed as follows,
cf. \cite{IsYo, TaYo}. There exists $\b_c \in (0,\8)$ such that
\bdmn 
\lim_{N \ra \8}
&Z_N^{1/N}& =\max \{ |\lm|\; ; \; \lm \in \s (P^+_V)\}\nn \\
& & \lef\{\barray{ll}
      = 1, & \mbox{if $\b  \le \b_c$ (delocalized phase)},\\
      >1, & \mbox{if $\b  > \b_c$ (localized phase)}. 
\earray \ri.
\label{wetting}
\edmn
Intuitively, if $N$ is large, the random walk under the measure $\m_N$ 
behaves as if there were no potentials in the delocalized phase,
while it is localized near $\6 \N^d$ in the localized phase.
In fact, for $d=1$, these intuitive pictures are justified mathematically
in \cite{IsYo}. In particular, it was shown there that, in the localized phase,
the measure $\m_N$ converges as $N \ra \8$ to the law of a positively
recurrent Markov chain. 

Technically, a crutial step in the proof of the limit theorem in \cite{IsYo} referred to above
is the existence of the absolute spectral gap of $P^+_V$ (cf. (\ref{asg}) below),
and it is here that the particularity of one dimension comes into play. 
Indeed, if $d=1$, the operator $P^+_V$ is a compact perturbation of $P^+$, from which the existence
of the absolute spectral gap follows via Weyl's
essential spectrum theorem.  
Unfortunately, $P^+_V-P^+$ is no longer a compact in higher dimensions.
The present article comes out as a partial progress in the effort to carry the results 
in \cite{IsYo} over to higher dimensions and therefore, to noncompact perturbation cases.

We now explain contents of this paper a little more in detail.
For simplicity, we consider the whole lattice $\zd$,
rather than its first quadrant $\N^d$.

Let $p:\zd \ra [0,\8)$ be a transition probability
  of a symmetric irreducible random walk. 
  We define $P:\ell^2(\Z^d) \ra \ell^2(\Z^d)$ by  
  \bdnl{P}
Pf(x)=\sum_{y \in \zd}p(x-y)f(y),\; \; f \in \ell^2(\zd),\; \; x \in \Z^d. 
\edn
We consider a perturbation $P_V:\ell^2(\zd) \ra \ell^2(\zd)$
of the operator $P$ of the form
  \bdnl{PV}
P_V= (1+V)P,
\edn
where $V:\Z^d \ra [0,\8)$ is a bounded function.
  Here, as usual, $V$ is regarded as the
  multiplication operator. In the context of probability theory,
  this type of perturbation is quite natural as we have already seen.
  
In this article, we investigate the structure of spectra $\s (P_V)$,
as well as its consequence on the long time behavior of the semigruop
(\ref{FKPV}). In this context, an important quantity is the
right edge
$$
r(P_V)=\max \s (P_V)
$$
of the spectra.
More precisely, we are interested in the following properties
\bdmn
r(P_V) & > & \sup\{\lm \; ; \; \lm \in \s (P_V),\; \lm \neq r(P_V)\}
\; \; \mbox{({\bf spectral gap})}, \label{sg}\\
r(P_V) & > & \sup\{|\lm| \; ; \; \lm \in \s (P_V),\; |\lm| \neq r(P_V)\}
\; \; \mbox{({\bf absolute spectral gap}).} \label{asg}
\edmn
As it turns out in the sequel, the following quantity plays an
important role in studying these properties. 
\bdnl{v0}
v_0=\inf_{n \in \N}\sup_{|x| \ge n}V(x). 
\edn
For example, as is discussed earlier in \cite{IsYo},
properties (\ref{sg}) and (\ref{asg}) are relatively easy to obtain
when $v_0=0$, i.e., $V$ decays at infinity. 
Indeed, the multiplication operator $V$ is compact in this case,  
and hence, by Weyl's essential spectrum theorem,
cf. \cite[p.358, Proposition 4.2 (e)]{Conway}, 
\cite[p.112, Theorem XIII.14]{RS78},
the set of essential spectra $\s_{\rm ess} (P_V)$
is the same as $\s_{\rm ess} (P)=\s(P)=[\ell (P),1]$
($-1 \le \ell (P)<1$). 
Thus, one immediately obtains (\ref{sg}) as soon as one knows that
$r(P_V)>1$. Then, it can be improved to (\ref{asg}) under
reasonable additional assumptions on the transition function $p$.

In this article, we are mainly interested in the case of
$v_0>0$, where the multiplication operator $V$ is no longer compact.
To compensate the lack of the decay of $V$ at infinity,
we will assume that the support of $V$ is sparse enough, to ensure
that the perturbation is not too large to control, see (\ref{sparse*1})
below for the precise formulation of the sparseness.
To the best of our knowledge, research in this direction
was initiated in \cite{Pearson}, where the Schr\"odinger
operator $-{d^2 \over dx^2}+V$ on the real line is discussed. 

Firstly, we characterize
in \Thm{ess} the set $\s_{\rm ess}(P_V)$ of essential
spectra of $P_V$.
Here, we adapt the method in \cite{Klaus} to the present setting.
\Thm{ess} has a corollary,
which tells us that the excess $\s_{\rm ess}(P_V) \bsh \s(P)$
is nonempty if, e.g., $d \le 2$ and $v_0>0$, cf. (\ref{v0}). 
For $d \ge 3$, we have the same conclusion if $v_0$ is sufficiently large. 
Thus, the spectral aspect of the operator $P_V$ is indeed different
from the compact perturbation case. 

Secondly, we prove in \Thm{decay} that that all the eigenfunctions which correspond to discrete spectra
  decay exponentially fast. 

Thirdly, we establish the spectral gap (\ref{sg}) in \Thm{gap}, and then,
improve it to the absolute spectral gap (\ref{asg}) in \Cor{gap}. 

Finally, as an application of the absolute spectral gap (\ref{asg})
obtained in \Cor{gap} and the exponential decay of the eigenfunctions (\Thm{decay}),
we prove a limit theorem for the
random walk in the potential $V$. More precisely, we consider a
Gibbs measures $\m_N$ ($N \in \N$) on $(\W,\cF)$ given by
\bdnl{Pn}
\m_N(d\w)={1 \over Z_N}E\lef[\prod^{N-1}_{j=0}(1+V(S_j)): d \w\ri], 
\edn 
where $Z_N$ is the normalizing constant. 
Then, we will prove that the measure $\m_N$ converges as $N \ra \8$ 
to a positively recurrent Markov chain on $\zd$, cf. \Thm{sup} below. 

The rest of this article is organized as follows. 
\tableofcontents 
\subsection{General notations}
\noindent $\bullet$
For a normed vector space $X$, we denote the totality
of bounded linear operators $T:X \ra X$ by $\cB (X)$.
For $T \in \cB (X)$, $\| T \|$, $\s (T)$ and $\rh (T)$ stands
respectively for its operator norm, the totality of its spectra,
and that of its resolvents. 
If $X$ is a Hilbert space and $T \in \cB (X)$ is self-adjoint,
we write  
\bdnl{r(T)ell(T)}
r(T)=\max \s (T),\; \; \; \ell (T)=\min \s (T).
\edn

\noindent $\bullet$
The Banach space $\ell^p (\Z^d)$ will be abbreviated by $\ell^p$. For a subset
$S \sub \zd$, $\ell^p (S)$ is identified with the totality of
$f \in \ell^p$ which vanish outside $S$.  

\noindent $\bullet$
For $T \in \cB (\ell^2)$,
we write its kernel by $T(x,y) \deq \lan \del_x, T \del_y\ran$
($x,y \in \Z^d$).

\noindent $\bullet$
For $u \in \ell^p$ and $u \in \ell^q$ ($p,q \in [1,\8]$, ${1 \over p}+{1 \over q}=1$)
let $\lan u,v \ran=\sum_{x \in \zd}u(x)v(x)^*$, where $c^*$ denotes the complex conjugate
of $c \in \C$. 
\subsection{Spectra of the operator $P_V$}
\label{sec:ess}
Let $p:\zd \ra [0,\8)$ be a transition probability
  of a random walk, i.e., $\sum_{x \in \zd}p(x)=1$. 
  Additionally, we assume that 
  \bdmn
  & & \mbox{({\bf symmetry}) $p(x)=p(-x)$ for all $x \in \zd$}; \label{symm}\\
  & & \mbox{({\bf finite range}) $p$ is supported on a finite set}; \label{f*range}\\
  & & \mbox{({\bf irreducibility}) for all $x \in \zd$, there exists
    $n \in \N$ such that $p_n(x)>0$,}
  \label{irred}
  \edmn
  where $p_n$ denotes the $n$-fold convolution. Then, we define
  $P \in \cB (\ell^2)$ and $P_V \in \cB (\ell^2)$ respectively
  by (\ref{P}) and (\ref{PV}).

  Let  $\widehat{p}(\tht)=\sum_{x \in \zd}p(x)\exp ({\bf i} \tht \cdot x)$,
$\tht \in [-\pi,\pi]^d$.  
Then, the set $\s(P)$ of the spectra of $P$ is the interval
$[\ell (P),1]$, where $\ell (P)= \min \widehat{p} \in [2p(0)-1,1)$.
For  $\lm \in \C \bsh \s(P)$, let $G_\lm \in \cB (\ell^2)$ be the
resolvent operator
\bdnl{Glm}
G_\lm = (\lm -P)^{-1}. 
\edn
We then define
\bdnl{g(lm)}
g_\lm(x) \deq \lm G_\lm(0,x), \; \; \lm \in \R \bsh \s(P),\; x \in \zd. 
\edn
In the sequel, we need to
deal with essential spectra of operators which are not self-adjoint.
To do so, we adopt the following definition via the theory of 
Fredholm operator. 
Let $X$ be a Banach space and $T \in \cB (X)$.
We say that $T$ is a Fredholm operator if
\bdnl{Fred}
\mbox{$\Ran T$ is closed, $\dim \Ker T <\8$, and $\dim \lef( X/\Ran T\ri)<\8$}.
\edn
We denote the totality of Fredholm operator by
${\rm F}(X)$.
We then define the set $\s_{\rm ess}(T) \sub \C$ of essential spectra of $T$ by
\bdnl{sig*ess}
\lm \in \s_{\rm ess}(T)
\; \;\Llra\; \;
\lm -T \not\in {\rm F}(X).
\edn
If $X$ is a Hilbert space and $T$ is self-adjoint, then, 
\bdnl{LSF4}
T \in {\rm F}(X)
\; \; \Llra \; \;
\dim \Ker T <\8
\;  \mbox{and}\;
0 \not\in \ov{\s(T) \bsh \{0\}}.
\edn
cf. \cite[p.359, Proposition 4.6]{Conway}.

In what follows, we will exploit
  the following characterization of the essential spectra of
  the multiplication operator $V$.
\bdnl{sig*ess*V}
v \in \s_{\rm ess}(V)
\; \;\Llra\; \;
\sharp \{ x \in \zd\; ; \; |V(x)-v| <\e\}=\8
\; \; \mbox{for all $\e>0$},
\edn
where $\sharp$ stands for the cardinality. 
Moreover, 
\bdnl{sig*ess*V2}
\max \s_{\rm ess}(V)=v_0, 
\edn
where $v_0$ is defined by (\ref{v0}). 
Let us consider an inner product and the associated norm. 
\bdnl{inner*prod*V}
\lan f, g \ran_V \deq \lan (1+V)^{-1}f, g \ran, 
\; \;\| f\|_V=\sqrt{\lan f, f \ran_V},
\; \; f,g \in \ell^2. 
\edn 
We note that the norms $\| \cdot \|$ and $\| \cdot \|_V$ are equivalent.
The Hilbert space $(\ell^2,\| \cdot \|_V)$ will be denoted by $\ell^2_V$. 
Then, the operator $P_V \in \cB (\ell^2_V)$ is self-adjoint.
As a consequence, all the spectra of $P_V$ are real numbers.

We first identify the set $\s_{\rm ess}(P_V)$ of the essential spectra. 
\Theorem{ess}
  Suppose that $V:\zd \ra [0,\8)$ is a bounded function such that 
\bdnl{sparse*1}
a_\e(x) \deq \sum_{y \in \zd \atop y \neq x}
  \sqrt{V(x)V(y)}\exp (-\e |x-y|)\st{|x| \ra \8}{\lra}0
\; \; \mbox{for any $\e>0$}. 
\edn
Then,
\bdnl{ess}
\s_{\rm ess}(P_V)=\s(P) \cup \Lm_V,
\edn
where 
\bdnl{Lm(v)1}
\Lm_V=\{ \lm \in \R \bsh \s(P)\; ;\;
\mbox{there exists $v \in \s_{\rm ess} (V) \bsh \{0\}$ such that $g_\lm (0)=1+v^{-1}$}\}.
\edn
\end{theorem}

\vs 
\noindent {\bf Remark} 
Condition (\ref{sparse*1}) is satisfied, not only
when $V(x)\st{|x| \ra \8}{\lra}0$, but also when
the support of $V$ is sufficiently sparse.
For example, the condition (\ref{sparse*1}) follows from (\ref{sparse*0}).
Indeed, suppose that $|x| \ge 2r$. Then, $|x-y| \ge r$ if $|y| \le r$.
Thus, letting ${\rm d}(r)$ denote the minimum on the left-hand side of (\ref{sparse*0}),
we have 
\bdnn
a_\e(x)
& \le &
\| V\|_\8\lef( \sum_{|y| \le r}+\sum_{|y| >r \atop y \neq x} \ri)\exp (-\e |x-y|)\\
& \le &
Cr^d\exp (-\e r)+C\exp \lef(-{\e {\rm d}(r) \over 2}\ri)
\sum_{y \in \zd \atop y \neq 0}\exp \lef(-{\e |y| \over 2}\ri)
\st{r \ra \8}{\lra}0. 
\ednn 

\vvs 
\Thm{ess} has the following corollary, which says that,
unlike the compact perturbation case, we may find essential
spectra of $P_V$ outside $\s(P)$. 
\Corollary{ess}
\bds
\item[a)]
$\s_{\rm ess}(P_V) \cap (1,\8) \neq \epty$, if and only if 
    \bdnl{cond*v}
\mbox{$v_0>0$ and $1+v_0^{-1} = g_{\lm_0} (0)$ for some $\lm_0 \in (1,\8)$},
\edn
where $v_0$ is defined by (\ref{v0}). Moreover, (\ref{cond*v}) implies
that $\lm_0=\max \s_{\rm ess}(P_V)$. 
\item[b)] $\s_{\rm ess}(P_V) \cap (-\8,\ell (P)) \neq \epty$,
if and only if
 \bdnl{cond*v2}
 \mbox{$\ell (P)<0$, $v_0>0$ and $1+v_0^{-1} = g_{\lm_0} (0)$ for some $\lm_0 \in (-\8, \ell (P))$}.
 \edn
 \eds 
\end{corollary}

\vs 
\noindent {\bf Remark} 
It follows from (\ref{g(lm)1}) that the condition (\ref{cond*v}) is satisfied
          whenever $\s_{\rm ess}(V) \bsh \{0\}\neq \epty$ if $d \le 2$. 
Similarly, we see from (\ref{g(lm)3}) that the condition (\ref{cond*v2}) is satisfied
          whenever $\s_{\rm ess}(V) \bsh \{0\}\neq \epty$ if $d \le 2$ and $\ell (P)<0$.  

          \vvs
          The next result (\Thm{decay}) deals with discrete spectra of $P_V$. More precisely, it
          states that the corresponding eigenfunctions decay exponentially fast.
          Together with the absolute spectral gap (\Cor{gap}), this theorem plays an important
          role in the proof of \Thm{sup} below. 
\Theorem{decay}
Suppose that the condition (\ref{sparse*1}) holds true, $\lm \in \s (P_V) \bsh \s_{\rm ess} (P_V)$
and that a function $\vp \in \ell^2$ satisfies $(\lm -P_V)\vp =0$. Then, there exist 
constants $\a,C \in (0,\8)$ such that
\bdnl{decay}
|\vp (x)| \le C\exp (-\a |x|),\; \; \mbox{for all $x \in \zd$}.
\edn
\end{theorem}

\vvs 
The next result deals with the right edge $r(P_V)$ of the spectra of $P_V$,
cf. (\ref{r(T)ell(T)}). 
\Theorem{gap}
In addition to conditions (\ref{sparse*1}) and (\ref{cond*v}), 
suppose that 
\bdnl{cond*v3}
V(\Z^d) \cap [v_0, \8) \neq \epty. 
  \edn
Then, there is a spectral gap at $r(P_V)$, i. e., (\ref{sg}) holds. Moreover, 
\bdmn
& & \mbox{there exists a strictly positive, normalized function $\vp \in \ell^2$ such that}\nn\\
& & \Ker (r(P_V)-P_V)=\C \vp.
\label{vp}
\edmn
\end{theorem}
  In the following corollary to \Thm{gap},
  we improve the spectral gap (\ref{sg}) to the absolute spectral gap (\ref{asg}).
  An operator 
  $T \in \cB (\ell^2)$ is said to be {\it bipartite} w.r.t. $J \in \{-1,1\}^{\zd}$ 
if $T(x,y)= 0$ for all $(x,y) \in (\zd)^2$ such that $J(x)J(y)=1$.
\Corollary{gap}
Suppose that (\ref{sg}) and (\ref{vp}) hold true, which is the case if
(\ref{sparse*1}), (\ref{cond*v}) and (\ref{cond*v3}) are satisfied. 
Then:
     \bds
     \item[a)]
       There is an absolute spectral gap at $r(P_V)$, i. e. (\ref{asg}) holds
       if and only if 
       \bdnl{cond*gap}
-r(P_V) <\ell (P_V)
\; \; \mbox{or $P$ is bipartite w.r.t. some $J \in \{-1,1\}^{\zd}$.} 
\edn
\item[b)]
  Assume (\ref{cond*gap}). 
Then, there exists $\e \in (0,1)$ such that
\bdnl{proj2}
\| r(P_V)^{-n}P^n_V (f -\Pi_V f)\|_V \le \e^n\|f -\Pi_V f\|_V
\; \; \;
\mbox{for all $f \in \ell^2$ and $n \in \N$},
\edn
where, with $\vp$ from (\ref{vp}),
\bdnl{proj}
\Pi_V f \deq 
\lef\{\barray{ll}
\lan f, \vp \ran_V \vp, 
& \mbox{if $-r(P_V) <\ell (P_V)$},\\
\lan f, \vp \ran_V \vp+\lan f, J\vp \ran_V J\vp,
& \mbox{if $P$ is bipartite w.r.t. $J \in \{-1,1\}^{\zd}$},
\earray \ri. 
\edn
\eds 
\end{corollary}
Finally, we provide a sufficient condition in terms of the
transition probability $p$ for the condition (\ref{cond*gap}).
\Proposition{cond*gap}
Let $A=\{x \in \zd\; ; \; \sum_{\a \in I}x_\a \in 2\Z\}$
with $\epty \neq I \sub \{1,\ldots,d\}$. 
Suppose that either $p$ vanishes on $A$, or
\bdnl{K*bsh*0}
p(0)>\sum_{x \in A \bsh \{0\}}p(x). 
\edn
Then, the condition (\ref{cond*gap}) is satisfied
for any nonnegative $V \in \ell^\8(\zd)$. More precisely, 
\bdnn 
& & \mbox{If $p$ vanishes on $A$, then, $P$ is bipartite w.r.t. $J={\bf 1}_A-{\bf 1}_{A^{\sf c}}$.}\\
& & \mbox{If (\ref{K*bsh*0}) holds, then,
  $-r(P_V)<\ell (P_V)$ for any nonnegative $V \in \ell^\8(\zd)$.} 
\ednn
  \end{proposition}
\subsection{Limit of the random walk in potential $V$}
\label{sec:mc}
We introduce a random walk $(S_n)_{n \in \N}$ on a probability space
$(\W, \cF, P)$ such that $S_0=0$, $P(S_1= \cdot )=p$. 
We define a probability measure $\m_N$ on $(\W, \cF)$ by (\ref{Pn}).

We assume (\ref{sparse*1}) and (\ref{cond*v}) in what follows.
We will have a positively recurrent Markov chain as the 
limit process of $\m_N$.
To describe the Markov chain obtained as the limit, 
we introduce an operator $P_{V,\vp} \in \cB (\ell^2)$ by 
\bdnl{Rxy}
P_{V,\vp}(x,y)
  =  r(P_V)^{-1}\vp (x)^{-1}P_V(x,y)\vp (y), 
\edn
where $\vp$ is from \Thm{gap}.
Then, $P_{V,\vp}$ is a transition probability for a Markov chain in $\zd$,   
which we denote by $(\{S_n\}_{n \in \N}, \{ \nu^x \}_{x \in \zd})$. 
It is easy to check that 
the Markov chain $\{ \nu^x \}_{x \in \zd}$ is 
reversible with respect to the probability measure $m$ on $\zd$ defined by
\bdnl{tlm(x)}
\lan f,m \ran=\frac{\lan f, \vp^2 \ran_V }{\lan {\bf 1}, \vp^2 \ran_V},
\; \; \; f \in \ell^\8. 
\edn 
In particular, 
$\{ \nu^x \}_{x \in \Z}$ is positively recurrent.
\Theorem{sup}
Assume the same hypothesis as \Cor{gap}. 
Then, there are constants $C=C(p,V)>0$ and
$\e=\e (p,V) \in (0,1)$ as follows; 
if $n \geq k +C$,  
$f: \Z^k \ra \R $ is polynomialy bounded and 
$F(\w)=f(S_1(\w), \ldots, S_k(\w))$, $\w \in \W$, then
\bdnl{muranu}
\lef|\int_\W Fd\mu_n-\int_\W Fd\nu \ri| \leq B(f) 
\e^{n-k},
\edn 
where $B(f)$ is a constant which depends 
only on $p$, $V$ and $f$. 
\end{theorem}
  Proof: Given \Thm{decay} and \Cor{gap}, the proof of \Thm{sup} is
  identical to that of \cite[Theorem 1.3]{IsYo},
  hence is omitted.
  \hfill $\Box$
          \subsection{An example in dimension one}
  \label{dim1}
We provide a simple example for $d=1$. 
We define $p:\Z \ra [0,1)$ such that $p(0)=q \in [0,1)$ $p(1)=p(-1)=(1-q)/2 \in (0,1/2]$. 
Then, for $\lm \in \R \bsh [2q-1,1]$, the function $g_\lm (x)$ is computed
explicitly. 
    \bdnl{green}
    g_\lm (x)
    ={\lm \over \sqrt{\del(\lm)}}
      \times \lef\{\barray{ll}
      \vp (\lm)^{|x|}, & \mbox{if $\lm >1$},\\
-\vp (\lm)^{-|x|}, & \mbox{if $\lm <2q-1$},
\earray \ri.
    \edn
where $\del(\lm)=(\lm -1)(\lm -(2q-1))>0$ and 
$$
\vp (\lm)={\lm -q -\sqrt{\del (\lm)} \over 1-q}
\in \lef\{\barray{ll}
(0,1), & \mbox{if $\lm >1$},\\
(-\8,-1), & \mbox{if $\lm <2q-1$}.
\earray \ri.
$$
\bdmn 
 & \bullet &\mbox{$g_\lm (0)$ is strictly decreasing in $\lm \in (1,\8)$, 
  $g_{\lm} (0)\st{\lm \searrow 1}{\lra}\8$, $g_\lm (0)\st{\lm \nearrow \8}{\lra}1$.}
\label{g(lm)1*1}\\
 & \bullet &\mbox{if $q=0$, then $g_\lm (0)$ is strictly increasing in $\lm \in (-\8,-1)$, 
  $g_{\lm} (0)\st{\lm \searrow -\8}{\lra}1$,}\nn \\
& & \mbox{$g_\lm (0)\st{\lm \nearrow -1}{\lra}\8$.}
\label{g(lm)2*1}\\
& \bullet & \mbox{if $0<q<1/2$,
  then, $g_\lm (0)$ is strictly decreasing in $\lm \in (-\8, {2q-1\over q})$,} \nn \\
& & \mbox{strictly increasing in $\lm \in ({2q-1\over q},2q-1)$,
  $g_{\lm} (0)\st{\lm \searrow -\8}{\lra}1$, $g_{2q-1\over q}(0)={\sqrt{1-2q} \over 1-q}$,}\nn \\
& & \mbox{$g_{2q-1\over 2q}(0)=1$, 
  $g_\lm (0)\st{\lm \nearrow 2q-1}{\lra}\8$},
\label{g(lm)3*1}\\
& \bullet & \mbox{ if $1/2 \le q <1$, then, $g_\lm (0)$ is strictly decreasing in
  $\lm \in (-\8,2q-1)$},\nn \\ 
& & \mbox{$g_{\lm} (0)\st{\lm \searrow -\8}{\lra}1$, $g_0 (0)=0$,
and if $1/2<q<1$, $g_\lm (0)\st{\lm \nearrow 2q-1}{\lra}-\8$.}  \label{g(lm)4*1}
\edmn
For $v>0$, we set
\bdnl{lmpm}
\lm_\pm (v) \deq c(v)\lef(q\pm \sqrt{q^2-(2q-1)c(v)^{-1}}\ri),
\; \; \mbox{where}\; \;
c(v) \deq {(v+1)^2 \over 2v+1}. 
\edn
Suppose that $V:\Z \ra [0,\8)$ is a bounded function which satisfies
  (\ref{sparse*1}).
  \bdnl{ess*1}
  \Lm_V=\lef\{\barray{ll}
\{\lm_\pm (v)\; ; \; v \in \s_{\rm ess}(V) \bsh \{0\}\},& \mbox{if $0 \le q <1/2$},\\
\{\lm_+ (v)\; ; \; v \in \s_{\rm ess}(V) \bsh \{0\}\}, & \mbox{if $1/2 \le q <1$}.
\earray \ri.
\edn
To see this, we 
  observe that for $v>0$ and $\lm \in \R \bsh [2q-1,1]$, 
  \bdnn
  g_\lm (0)=1+v^{-1}&\Lra&{\lm^2 \over \del (\lm)}=(1+v^{-1})^2 \\
  &\Llra& \lm^2-2qc(v)\lm +(2q-1)c(v)=0 \\
  &\Llra& 
\lm =\lm_\pm (v).
\ednn
Moreover, 
\bdnl{lmpm2}
\lm_{-}(v)<2q-1<1<
\lm_{+}(v)=\lm_{-}(v)+2c(v)\sqrt{q^2-(2q-1)c(v)^{-1}}.
\edn
Taking (\ref{g(lm)1*1})--(\ref{lmpm2}) into account
(In particular, if $0<q<1/2$, then, ${2q-1 \over 2q} <\lm_{-}(v)<2q-1$ and hence
$g_{\lm}(0)>1$ at $\lm=\lm_{-}(v)$ by (\ref{g(lm)3*1})), we have 
$$
\lm \in \R \bsh [2q-1,1],\; g_\lm (0)=1+v^{-1}
\; \Llra \;
\lm =\lef\{\barray{ll}
\lm_\pm (v),& \mbox{if $0 \le q <1/2$},\\
\lm_+ (v), & \mbox{if $1/2 \le q <1$}.
\earray \ri.
$$
This, together with \Thm{ess}, implies (\ref{ess*1}).
  \SSC{Proof of \Thm{ess} and \Cor{ess}}
\subsection{Outline}
\noindent {\bf Proof of \Thm{ess}}
{\it Step1} We prove that $\s_{\rm ess}(P_V) \bsh \s(P)=\Lm_V$. 
Here, we adapt the method in \cite{Klaus} to the present setting.
Let $\lm \in \R \bsh \s(P)$. 
Then,
We introduce the modified Birman-Schwinger
operator $G_{V,\lm} \in \cB (\ell^2)$ by
\bdnl{BS}
G_{V, \lm}= V^{1/2}(\lm G_\lm -1)V^{1/2}. 
\edn
We first show in \Lem{BS} that
\bdnl{lem*BS}
\lm \in \s_{\rm ess}(P_V)
\; \; \Llra\; \;  
1 \in \s_{\rm ess}(G_{V,\lm}).
\edn
We then prove in \Lem{Weyl} that
\bdnl{lem*Weyl}
G_{V,\lm}=(g_\lm (0)-1)V+H_{V,\lm},
\edn 
where $H_{V,\lm}$ is a compact operator. In fact, this
is where the condition (\ref{sparse*1}) is used. 
Then, we see from (\ref{lem*Weyl}) that
\bdnl{lem*Weyl*2}
\s_{\rm ess}(G_{V,\lm})
=(g_\lm (0)-1)\s_{\rm ess}(V), 
\edn 
via Weyl's essential spectrum theorem,
cf. \cite[p.358, Proposition 4.2 (e)]{Conway}, 
\cite[p.112, Theorem XIII.14]{RS78}.
By (\ref{lem*Weyl*2}), 
\bdnl{lem*Weyl*3}
1 \in \s_{\rm ess}(G_{V,\lm})
\; \; \Llra\; \;
\exists v \in \s_{\rm ess}(V) \bsh \{0\},\;
g_\lm (0)-1=v^{-1} 
\; \; \Llra\; \;
\lm \in \Lm_V.
\edn 
Thus, we obtain (\ref{ess}) from (\ref{lem*BS}) and (\ref{lem*Weyl*3}). \\
{\it Step2} We prove that $\s(P) \sub \s_{\rm ess}(P_V)$. This step
is the subject of \Lem{sub*ess} below. 
\hfill $\Box$

\vvs
\noindent {\bf Proof of \Cor{ess}}:
a) Suppose that the condition (\ref{cond*v}) holds. 
Then, $\lm_0 \in \Lm_V \cap (1,\8)$.
Since $\s_{\rm ess}(V) \cap (1,\8)=\Lm_V \cap (1,\8)$ by \Thm{ess}, 
it follows that $\s_{\rm ess}(V) \cap (1,\8) \neq \epty$.
Moreover, by (\ref{g(lm)1}), $g_{\lm}(0)$ is strictly decreasing in $\lm \in (1,\8)$. 
Therefore. (\ref{cond*v}) implies that $\lm_0=\max \s_{\rm ess}(P_V)$. \\
Suppose on the other hand that $\s_{\rm ess}(V) \cap (1,\8) \neq \epty$.
Then, since $\s_{\rm ess}(V) \cap (1,\8)=\Lm_V \cap (1,\8)$ by \Thm{ess}, 
there exist $\lm \in (1,\8)$ and $v \in \s_{\rm ess}(V) \bsh \{0\}$ such that $g_{\lm}(0)=1+v^{-1}$.
Then, $v_0 \ge v >0$. Moreover, since $0 <v_0^{-1} \le v^{-1}$ and by (\ref{g(lm)1}),  
$g_\lm (0) \st{\lm \ra \8}{\lra}1$, 
it follows from the mean value theorem that
there exists $\lm_0 \in (1,\8)$ such that $g_{\lm_0}(0)=1+v_0^{-1}$. \\
b) The proof is similar as above. 
\hfill $\Box$
\subsection{Lemmas}
Properties of the function $\lm \mapsto g_\lm (0)$, which we need in this article are summarized
in the following 
\Lemma{g(lm)}
\bdmn 
 & \bullet &\mbox{$g_\lm (0)$ is strictly decreasing in $\lm \in (1,\8)$, 
  $g_\lm (0)\st{\lm \ra \8}{\lra}1$.}\nn \\
& &\mbox{Moreover, $g_{1+} (0)\lef\{\barray{ll}
      =\8, & \mbox{if $d \le 2$},\\
      \in (0,\8), & \mbox{if $d \ge 3$}.
      \earray \ri.$}\label{g(lm)1}\\
& \bullet & \mbox{if $\ell (P)<0$, then, $g_\lm (0)>{|\lm| \over |\lm|+p(0)}$
  for $\lm \in (-\8,\ell (P))$.} \nn \\
& & \mbox{$g_\lm (0)\st{\lm \ra -\8}{\lra}1$, Moreover, $g_{\ell (P)-}(0)=\8$ if $d \le 2$},\label{g(lm)3}\\
& \bullet & \mbox{ if $\ell (P)\ge 0$, then, $g_\lm (0)$ is strictly decreasing in
  $\lm \in (-\8,\ell (P))$},\nn \\ 
  & & \mbox{$g_\lm (0) \in\lef\{\barray{ll}
       (0,1), & \mbox{for $\lm \in (-\8,0)$},\\
      (-\8,0), & \mbox{for $\lm \in (0,\ell (P))$}.
      \earray \ri.$}  \label{g(lm)4}
  \edmn
  \end{lemma}
  Proof: The behavior (\ref{g(lm)1}) of $g_\lm (0)$ for $\lm >1$ is 
  well-known in the context of the random walk. Thus, we omit the proofs.
  
  On the other hand, the following integral formula is well-known. 
\bdnl{lmGlm}
g_\lm (0)
={\lm \over (2\pi)^{d}}\int_{[-\pi,\pi]^d}{ d\tht \over \lm -\widehat{p}(\tht)},
\; \; \; \lm \in \R \bsh \s(P).
\edn
We derive properties (\ref{g(lm)1})--(\ref{g(lm)4}) from this formula.
We take (\ref{g(lm)3}) for example. Since $\lm$ is negative, 
\bdnn
g_\lm (0)
&=&{|\lm| \over (2\pi)^{d}}\int_{[-\pi,\pi]^d}{ d\tht \over |\lm| +\widehat{p}(\tht)}\\
&>&{|\lm| \over |\lm| +(2\pi)^{-d}\int_{[-\pi,\pi]^d}\widehat{p}(\tht)d\tht}
={|\lm| \over |\lm|+p(0)},
\ednn
where we have used Jensen inequality to the
convex function $x \mapsto 1/x$ ($x>0$)
to proceed from the
first line to the second. 

To show that $g_{\ell (P)-}(0)=\8$ if $d \le 2$, we
take $\tht_0 \in [-\pi,\pi]^d$ such that $\ell (P)=\widehat{p}(\tht_0)$.
Since ${\6 \widehat{p} \over \6 \tht_\a} (\tht_0)=0$ for all
$\a=1,\ldots,d$, we have  
\bdnn 
\widehat{p}(\tht_0+\tht)-\ell (P)
&=&
\widehat{p}(\tht_0+\tht)-\widehat{p}(\tht_0) \\
&=&{1 \over 2}\sum_{\a,\b=1}^d\tht_\a \tht_\b\int^1_0(1-t)
{\6^2 \widehat{p} \over \6 \tht_\a\6 \tht_\b} (\tht_0 +t \tht)dt \\
& \le & C|\tht|^2
\; \; \;\mbox{for $\tht \in [-\pi,\pi]^d$}. 
\ednn
By applying the monotone convergence theorem to the integral
$\int_{[-\pi,\pi]^d}{ d\tht \over \widehat{p}(\tht)-\lm}$
as $\lm \nearrow \ell (P)$, we have
\bdnn
g_{\ell (P)-}(0)
&=& {|\ell (P)| \over (2\pi)^{d}}\int_{[-\pi,\pi]^d}{ d\tht \over \widehat{p}(\tht)-\ell (P)}
= {|\ell (P)|\over (2\pi)^{d}}\int_{[-\pi,\pi]^d}{ d\tht \over \widehat{p}(\tht_0+\tht)-\ell (P)}\\
& \ge &
{|\ell (P)| \over (2\pi)^{d}C}\int_{[-\pi,\pi]^d}{ d\tht \over |\tht|^2}=\8.
\ednn 
\hfill $\Box$

\vvs
\noindent {\bf Remark}
If $\ell (P)<0$, $g_\lm (0)$ is not necessarily monotone in $\lm \in (-\8,\ell (P))$.
See (\ref{g(lm)3*1}) for a counterexample. 

\Lemma{AB}
Let $A,B \in \cB (X)$ on a Banach space $X$.
Then,
\bds
\item[a)]
  $1 \in \s_{\rm ess}(BA)$ if and only if $1 \in \s_{\rm ess}(AB)$.
\item[b)]
  $1 \in \rh (BA)$ if and only if $1 \in \rh (AB)$.
  Moreover, if these conditions hold true, then   
  \bdnl{AB}
(1-BA)^{-1}=1+B(1-AB)^{-1}A. 
  \edn 
  \eds 
\end{lemma}
  Proof: a) Since the roles of $A$ and $B$ are exchangeable,
  it is enough to verify only the "if'' part, which is
  equivalently stated as $1-AB \in {\rm F}(X)$ $\Ra$ $1-BA \in {\rm F}(X)$. 
To prove this, we use Atkinson' s theorem which says 
the following for $T \in \cB (X)$, cf. \cite[p.161, Theorem 4.46]{AbAl}. 
\bds 
\item[A1)] Suppose that $T \in {\rm F}(X)$ with
  $n_1=\dim \Ker T$ and $n_2=\dim (X /\Ran T)$.
  Then, there exist $S,K_1,K_2 \in \cB(X)$ with 
  ${\rm rank} K_j=n_j$ ($j=1,2$) such that
  \bdnl{Atkinson}
  ST=1+K_1\; \; \mbox{and}\; \;  TS=1+K_2. 
  \edn 
\item[A2)] Conversely, suppose that there exist $S,K_1,K_2 \in \cB(X)$
  of which $K_1$ and $K_2$ are compact such that (\ref{Atkinson}) holds true.
  Then, $T \in {\rm F} (X)$.
  \eds  
  Comming back to the proof of the lemma, suppose that
  $T \deq 1 -AB \in {\rm F} (X)$. Then, by Atkinson's theorem, 
there exist $S,K_1,K_2 \in \cB(X)$ as are stated in A1) above. 
Then, $1+BSA \in \cB(X)$, ${\rm rank}BK_jA \le n_j$ ($j=1,2$) and 
\bdnl{AB2}
(1+BSA)(1-BA)=1+BK_1A,
\; \; 
(1-BA)(1+BSA)=1+BK_2A. 
\edn 
Therefore, $1 -BA \in {\rm F} (X)$, by A2) above. 
\\ b) The equivalence of $1 \in \rh (BA)$ and $1 \in \rh (AB)$
can be regarded as a special case of part a), where $K_1=K_2=0$
in the proof above. Moreover, suppose that $1 \in \rh(BA)$,
or equivalently, $1 \in \rh(AB)$. Then, $S=(1-AB)^{-1}$ in
(\ref{AB2}). Thus, the equality (\ref{AB}) follows from (\ref{AB2}). 
\hfill $\Box$ 
\Lemma{BS}
Let $\lm \in \R \bsh \s(P)$, and $G_{V,\lm}$ be the operator defined by (\ref{BS}).
Then,
\bds
\item[a)]
  $\lm \in \s_{\rm ess}(P_V)$ if and only if $1 \in \s_{\rm ess}(G_{V,\lm})$.
\item[b)]
  $\lm \in \s(P_V)$ if and only if $1 \in \s(G_{V,\lm})$.
    Moreover, if $\lm \in \rh (P_V)$, or equivalently, $1 \in \rh (G_{V,\lm})$, then, 
  \bdnl{eq:BS}
(\lm - P_V)^{-1}=G_\lm+G_\lm V^{1/2}(1-G_{V,\lm})^{-1}V^{1/2}PG_\lm. 
  \edn 
\eds
\end{lemma}
  Proof: a) We need to verify the equivalence.  
  $\lm -P_V \in F(\ell^2)$ $\Llra$ $1 -G_{V,\lm} \in F(\ell^2)$.
  We decompose this task into the following two steps.   
  \bdmn
  \lm -P_V \in F(\ell^2)\; \; &\Llra& \; \;
1-G_\lm PV  \in F(\ell^2),
\label{BS*1}\\
1-G_\lm PV \in F(\ell^2)\; \; &\Llra& \; \;
1 -G_{V,\lm} \in F(\ell^2).
\label{BS*2}
\edmn
Note that
\bdnl{algebra}
\lm -P_V =\lm -P -VP =(\lm -P)(1-G_\lm VP).
\edn
or equivalently,
\bdnl{algebra2}
1-G_\lm VP=G_\lm (\lm -P_V).
\edn
Note that $\lm -P, G_\lm \in F(\ell^2)$. Recall also that $F(\ell^2)$
is closed under the composition cf. \cite[p.158, Theorem 4.43]{AbAl}.
Then, (\ref{BS*1}) follows from
(\ref{algebra}) and (\ref{algebra2}). 

The equivalence (\ref{BS*2}) follows from \Lem{AB} a), since
for $A \deq G_\lm V^{1/2}$ and $B \deq V^{1/2}P$, 
\bdnn
AB &=& G_\lm V^{1/2}V^{1/2}P=G_\lm VP,\\
BA &=& V^{1/2}PG_\lm V^{1/2}=V^{1/2}(\lm G_\lm-1)V^{1/2}=G_{V,\lm}.
\ednn
b) The proof is similar as above. We use \Lem{AB} b) instead of \Lem{AB} a).
Suppose that $\lm \in \rh(P_V)$. Then, it follows from (\ref{algebra2}) that
$$
(\lm - P_V)^{-1}=(1-G_\lm VP)^{-1}G_\lm. 
$$
On the other hand, by plugging $A = V^{1/2}P$ and $B =  G_\lm V^{1/2}$ into
(\ref{AB}),
\bdnn
(1-G_\lm VP)^{-1}
&=& 1+G_\lm V^{1/2}(1-V^{1/2}PG_\lm V^{1/2})^{-1}V^{1/2}P\\
&=& 1+G_\lm V^{1/2}(1-G_{V,\lm})^{-1}V^{1/2}P.
\ednn
Combining these, we obtain (\ref{eq:BS}). 
\hfill $\Box$
\Lemma{Weyl}
  With $g_\lm (0)$, $G_{V,\lm}$ defined respectively by (\ref{g(lm)}) and (\ref{BS}),
  the following operator $H_{V,\lm} \in \cB (\ell^2)$ is compact.
  $$
H_{V,\lm} \deq G_{V,\lm}-(g_\lm (0)-1)V.
$$
\end{lemma}

  Proof:
We decompose $G_{V,\lm}$ into diagonal and off diagonal components as follows.
  $$
  G_{V,\lm}=V^{1/2}(\lm G_\lm -1)V^{1/2}= (g_\lm (0)-1)V+H_{V,\lm},
  $$
  where the operator $H_{V,\lm}$ is given by the kernel
  $$
  H_{V,\lm}(x,y)
  =\lm \sqrt{V(x)V(y)}G_\lm(x,y){\bf 1}_{x \neq y}.
$$
  We will show that
  \bdnl{K-KN}
\|H_{V,\lm}-H_{V,\lm}^{(N)}\| \st{N \ra \8}{\lra}0,
  \edn 
  where $K_N$ is a finite rank operator  defined by the kernel
  $$
  H_{V,\lm}^{(N)}(x,y)
  =\lm \sqrt{V(x)V(y)}G_\lm(x,y){\bf 1}_{x \neq y, |x| \le N}.
  $$
  Thus, (\ref{K-KN}) shows that $H_{V,\lm}$ is a compact operator.

  By a standard estimate, 
  $$
  \|H_{V,\lm}-H_{V,\lm}^{(N)}\| \le |\lm| \sqrt{A_NB_N},
  $$
  where
  $$
  A_N
  =\sup_{|x| \ge N}\sum_{y \in \zd \atop x \neq y}
  \sqrt{V(x)V(y)}|G_\lm(x,y)|, 
  \; \; \; 
  B_N
  =\sup_{x \in \zd}\sum_{y \in \zd \atop y \neq x, |y| \ge N}
  \sqrt{V(x)V(y)}|G_\lm(x,y)|.
  $$
  By (\ref{f*range}), there exist
  $C=C(\lm) \in (0,\8)$ and $\e=\e (\lm) \in (0,\8)$ such that
  \bdnl{green*dec}
|G_\lm(x,y)|\le C\exp (-\e |x-y|). 
\edn
Thus, (\ref{K-KN}) follows from (\ref{sparse*1}) and 
\bdnl{bN*ra*0}
\sup_{x \in \zd}b_{N,\e}(x)\st{N \ra \8}{\lra}0,
\; \; \mbox{where}\; \; 
b_{N,\e}(x)=\sum_{y \in \zd \atop y \neq x, |y| > N}
  \sqrt{V(x)V(y)}\exp (-\e |x-y|). 
\edn 
Therefore, it is enough to show that (\ref{sparse*1}) implies
(\ref{bN*ra*0}). Moreover, it follows immediately from (\ref{sparse*1})
that $\sup_{|x| >\lfloor N /2 \rfloor}b_{N,\e}(x)\st{N \ra \8}{\lra}0$. 
Hence, it is enough to verify that
\bdnl{bN*ra*1}
\sup_{|x| \le \lfloor N /2 \rfloor}b_{N,\e}(x)\st{N \ra \8}{\lra}0.
\edn 
If $|x| \le \lfloor N /2 \rfloor$ and $|y| \ge N$, then $|y-x| \ge N/2$.
Thus,
\bdnn
\sup_{|x| \le \lfloor N /2 \rfloor}b_{N,\e}(x)
& \le & \exp (-\e N /2)\|V\|_\8
\sup_{|x| \le \lfloor N /2 \rfloor}\sum_{y \in \zd}\exp (-\e |x-y|/2)\\
& = & \exp (-\e N /2)
\|V\|_\8\sum_{y \in \zd}\exp (-\e |y|/2)
\st{N \ra \8}{\lra}0,
\ednn
which proves (\ref{bN*ra*1}).
\hfill $\Box$ 

\vvs
In what follows, we use the following notation.
For $x=(x_\a)_{\a=1}^d \in \zd$,
\bdnl{def*||8}
|x|_\8 =\max_{1 \le \a \le d}|x_\a|.
\edn
For $c \in \zd$ and a positive integer $\ell$,
\bdnl{def*cube}
Q(c,\ell)=\{x \in \zd\; ; \; |x-c|_\8 \le \ell\}.
\edn
Suppose that the condition (\ref{sparse*1}) holds true
and that $v_0>0$. Then the following lemma shows that
there are infinitely many disjoint cubes
$Q(c,\ell)$ in which $V$ takes the value close to $v_0$ at the
center $c$, while the value $V(x)$ for the other points  
of the cube are close to zero.  
\Lemma{high}
Suppose that the condition (\ref{sparse*1}) holds true
and that $v_0>0$.
Then, for any $L, \ell \in (0,\8)$ and
$\e \in (0,1)$, there exists $c \in \zd$ such that
\bdmn
Q(0,L) \cap Q(c,\ell)&=&\epty,
\label{high1}\\
V(c) & > & (1-\e)v_0,
\label{high2}\\
\sum_{x \in Q(c,\ell)\bsh \{c\}}V(x) & < & \e.
\label{high3}
\edmn
\end{lemma}
Proof: Take $\del \in (0,1)$ such that $\del^2\exp (2\ell )<\e(1-\e)v_0$. 
By assumption $v_0>0$, the set
$$
H_\e =\{ x \in \zd \; ; \; V(x) >(1-\e)v_0\}
$$
is unbounded. Therefore, by (\ref{sparse*1}), we can find $c \in H_\e$
such that
\bdnl{high4}
|c|>L+\ell,\; \; \mbox{and}\; \; 
A(c) \deq \sum_{x \in \zd \atop x \neq c}
  \sqrt{V(c)V(x)}\exp (-|c-x|_\8)<\del.
  \edn
  The first inequality of
  (\ref{high4}) is equivalent to (\ref{high1}),
  while the second inequality implies (\ref{high3}) as follows.
  \bdnn
\lefteqn{  \sum_{x \in Q(c,\ell)\bsh \{c\}}V(x)} \\
   &\le &
  \lef(  \sum_{x \in Q(c,\ell)\bsh \{c\}}\sqrt{V(x)}\ri)^2
   \st{\scriptsize \mbox{(\ref{high2})}}{\le} 
  A(c)^2{\exp (2\ell) \over (1-\e)v_0}
  \st{\scriptsize \mbox{(\ref{high4})}}{<}
  \del^2{\exp (2\ell) \over (1-\e)v_0}<\e.
  \ednn
  \hfill $\Box$ 
\Lemma{sub*ess}
Suppose that the condition (\ref{sparse*1}) holds true.
Then, $\s(P) \sub \s_{\rm ess}(P_V)$. 
\end{lemma}
Proof:
Let $\lm \in \s(P)$ be arbitrary. We will prove that $\lm -P_V \not\in {\rm F}(\ell^2)$,
or equivalently, $\lm -P_V^* \not\in {\rm F}(\ell^2)$, where $P_V^*$ is the adjoint
operator of $P_V$ on $\ell^2$, therefore, $P_V^*=P+PV$.

For this purpose, we will use the following criterion
for an operator $T$ on a Hilbert space $X$ not to be a Fredholm operator.
$T \not\in {\rm F}(X)$ if there exists a normalized
sequence $\{u_n\} \sub X$  such that 
\bdnl{LSF2}
\mbox{$u_n \st{n \ra \8}{\lra}0$ weakly and $Tu_n \st{n \ra \8}{\lra}0$ strongly,}
\edn
The converse is also true if $T$ is self-adjoint, cf. \cite[p. 350, Theorem 2.3]{Conway}. 
A sequence with above property is called a {\it Weyl sequence}.

We will construct a Weyl sequence $u_n \in \ell^2$
for $\lm -P_V^*$, that is, $u_n \in \ell^2$ is normalized,
weakly convergent to zero, and 
\bdnl{sub*ess*Weyl}
\|(\lm -P-PV)u_n\|\st{n \ra \8}{\lra}0.
\edn
{\it Step1}: 
Let $r>0$ be an integer such that ${\rm supp}[p] \sub Q(0,r)$,
We first show that, there exist $c_n \in \zd$, $n \in \N$ such that
\bdmn
Q(c_m,m+r) \cap Q(c_n,n+r)&=&\epty
\; \; \; \; \mbox{if $m \neq n$},
\label{cube1}\\
\max_{x \in Q(c_n,n+r)}V(x) &<& (n+2)^{-1}
\; \; \mbox{for all $n \in \N$}.
\label{cube2}
\edmn
Suppose that $v_0=0$. Then, for any $\e>0$, there exists $R>0$ such that
$V(x) < \e$ if $|x| >R$. Thus, it is easy to
find such $c_n$ by an inductive procedure.

Suppose on the contrary that $v_0>0$. We then proceed inductively
with the help of \Lem{high} as follows. 
We start by taking $L=1$, $\ell=3r$ and $\e=1/2$ in \Lem{high},
so that we can find $b_0 \in \zd$ such that
$$
\max_{x \in Q(b_0,3r)\bsh \{b_0\}}V(x) < 1/2.
$$
Thus, by choosing $c_0$ such that $Q(c_0,r) \sub Q(b_0,3r)\bsh \{b_0\}$,
we obtain (\ref{cube2}) for $n=0$.
Next, we take $L=|c_0|+1+r$, $\ell=3+3r$ and $\e=1/3$ in \Lem{high},
so that we can find $b_1 \in \zd$ such that
$$
Q(c_0, r) \cap Q(b_1, 3+3r) =\epty,
\; \; \mbox{and}\; \;
\max_{x \in Q(b_1,3+3r)\bsh \{b_1\}}V(x) < 1/3.
$$ 
Thus, by choosing $c_1$ such that $Q(c_1,1+r) \sub Q(b_1,3+3r)\bsh \{b_1\}$,
we obtain (\ref{cube2}) for $n=1$. By repeating this procedure, we obtain
$c_n$, $n \in \N$ as desired. 

\vs
{\it Step2}: We now construct the Weyl sequence $u_n$, $n \in \N$. 
Let $\tht \in [-\pi,\pi]^d$ be such that $\lm=\widehat{p}(\tht)$ and
set $e_\lm (x)=\exp (\bfi \tht \cdot x)$. 
Then by Step1, there exist $c_n \in \zd$, $n \in \N$ which satisfy 
(\ref{cube1}) and (\ref{cube2}). We set
\bdnl{sub*ess*u}
\vp_n=e_\lm {\bf 1}_{Q(c_n,n+r)}
\; \; \mbox{and}\; \;
u_n=\vp_n/\|\vp_n\|.
\edn
We will show that $u_n$, $n \in \N$ is the Weyl sequence which we look for. 
By (\ref{cube1}), $u_n$, $n \in \N$ are orthonormal.
Moreover, we see from (\ref{cube2}) that
$$
0 \le V \le (n+2)^{-1}{\bf 1}_{Q(c_n,n+r)},
$$
and hence $\|Vu_n\| \le (n+2)^{-1}\st{n \ra \8}{\lra}0$. 
Thus, it only remains to verify that
\bdnl{sub*ess*Wey2}
\|(\lm -P)u_n\|\st{n \ra \8}{\lra}0.
\edn
To see this, we observe that
\bdnl{sub*ess*ev}
(\lm -P)e_\lm= 0.
\edn
Note also that for $f,g: \zd \ra \C$ and $x \in \rd$, 
$$
\mbox{$f=g$ on $Q(x,r)$}\; \; \Lra \; \; Pf(x)=Pg(x). 
$$
This, together with (\ref{sub*ess*u}) and (\ref{sub*ess*ev}), implies that
\bdnl{sub*ess*ev2}
\mbox{$(\lm-P) \vp_n=0$ on $Q(c_n,n)$}.
\edn
If we set $h_n \deq {\bf 1}_{Q(c_n,n+r)\bsh Q(c_n,n)}$, then
$$
\|h_n\|^2=\sharp (Q(c_n,n+r)\bsh Q(c_n,n))=(2n+2r+1)^d-(2n+1)^d=O(n^{d-1}),
$$
and hence
\bdnl{sub*ess*Wey3}
\|(\lm -P)\vp_n\|
\st{\scriptsize \mbox{(\ref{sub*ess*ev2})}}{=}
\|h_n(\lm -P)\vp_n\|
\le
\|h_n\|\|(\lm -P)\vp_n\|_\8=O(n^{(d-1)/2}).
\edn 
This implies (\ref{sub*ess*Wey2}), since 
$\|\vp_n\|=\|{\bf 1}_{Q(c_n,n+r)}\|=(2n+2r+1)^{d/2}$. 
\hfill $\Box$ 
\SSC{Proof of \Thm{decay}}
\subsection{Outline}\label{decay:outline}
For $\a>0$,
we denote by $\ell^{\8,\a}$ the Banach space of 
exponentially decaying functions $u:\zd \ra \C$ with the exponent at least $\a$,
more precisely, the functions which satisfy  
\bdnl{ell8e}
\|u\|_{\8,\a} \deq \sup_{x \in \zd}\exp (\a |x|)|u(x)|<\8.
\edn
We also recall the estimate (\ref{green*dec})
on the exponential decay of the kernel $G_\lm (x,y)$. 

Suppose that the condition (\ref{sparse*1}) holds true.
Let $\lm \in \C \bsh \s_{\rm ess} (P_V)$ and $\a \in (0,\e)$, where
$\e>0$ is from (\ref{green*dec}). For $K \sub \zd$, we define
\bdnl{VW}
V_K=V{\bf 1}_{\zd \bsh K}.
\edn
In the sequel, we write $K \Subset \zd$, when $K$ is a finite subset of $\zd$. 
We will prove in \Lem{P*W} below that there exists $K \Subset \zd$ such that
$\lm \in \rh (P_{V_K})$ and $(\lm -P_{V_K})^{-1} \in \cB (\ell^{\8,\a})$.  

Suppose additionally that $\lm \in \s (P_V)$ 
and that a function $\vp \in \ell^2$ satisfies $(\lm -P_V)\vp =0$.
Then, with the set $K$ from \Lem{P*W}, we rewrite
$(\lm -P_V)\vp =0$ as 
$$
(\lm -P_{V_K})\vp ={\bf 1}_KVP\vp.
$$
Since $\lm \in \rh (P_{V_K})$ by \Lem{P*W}, it follows from the above display that 
\bdnl{decay2}
\vp =(\lm -P_{V_K})^{-1}{\bf 1}_KVP\vp. 
\edn 
Since the function ${\bf 1}_KVP\vp$ is supported on the finite set $K$
and $(\lm -P_{V_K})^{-1} \in \cB (\ell^{\8,\a})$ by \Lem{P*W},
we obtain (\ref{decay}) from (\ref{decay2}).
\hfill $\Box$ 
\subsection{Lemmas}
Let $\lm \in \C \bsh \s_{\rm ess} (P_V)$ and $\a \in (0,\e)$, where
$\e>0$ is from (\ref{green*dec}). 
As is discussed in section \ref{decay:outline}, it is enough to prove 
that there exists $K \Subset \zd$ such that
$\lm \in \rh (P_{V_K})$ and $(\lm -P_{V_K})^{-1} \in \cB (\ell^{\8,\a})$
(\Lem{P*W}).
We will implement this by dealing with the modified Birman-Schwinger
operator with the potential $V_K$:
$$
G_{V_K,\lm}=V_K^{1/2}(\lm G_\lm -1)V_K^{1/2},
$$
cf. (\ref{BS}). 
As in the proof of \Lem{Weyl},
we decompose $G_{V_K,\lm}$ into diagonal and off diagonal components as follows.
  \bdnl{G*W*decomp}
  G_{V_K,\lm}= \gm V_K+H_{V_K,\lm},
  \; \; \mbox{where}\; \;\gm =g_\lm (0)-1. 
  \edn
  We first look at the diagonal component $\gm V_K$ of the above decomposition.
  \Lemma{1-gm*W}
Suppose that $\lm \not\in \s (P) \cup \Lm_V$. 
Then, there exists $K_0 \Subset \zd$ 
such that if $K_0 \subset K \Subset \zd$, 
\bdnl{1-gm*W}
\e_0 \deq \inf_{x \in \zd}|1-\gm V_K(x)|>0. 
\edn
\end{lemma}

\noindent Proof: 
It follows from the assumption $\lm \not\in \Lm_V$ that either 
$$
\mbox{{\bf i)} $\gm =0$}
\; \;\mbox{or}\; \;
\mbox{{\bf ii)} $\gm \neq 0$ and $\gm^{-1} \not\in \s_{\rm ess}(V)$}.
$$
If $\gm=0$, then, (\ref{1-gm*W}) is clearly true with $K =\epty$ and $\e_0=1$.
Suppose ii) above. 
Then, there exists $\del >0$ such that 
$$
K_0 \deq \{ x \in \zd \; ; \; |1-\gm V(x)|<\gm \del \}
=\{ x \in \zd \; ; \; |\gm^{-1}-V(x)|<\del \}
\Subset \zd.
$$
Let $K_0 \sub K \Subset \zd$. Then, it is clear that
$$
|1-\gm V_K(x)|=|1-\gm V(x)| \ge \gm \del \; \; \mbox{for all $x \in \zd \bsh K$}, 
$$
which implies (\ref{1-gm*W}). 
\hfill $\Box$

\vvs
The following lemma deals with the off diagonal part $H_{V_K,\lm}$ 
of the operator $G_{V_K,\lm}$, cf. (\ref{G*W*decomp}) which is given
by the kernel.
  \bdnl{H*W*ker}
  H_{V_K,\lm}(x,y)
  =\lm \sqrt{V_K(x)V_K(y)}G_\lm(x,y){\bf 1}_{x \neq y}.
  \edn 
  \Lemma{H*W}
  Suppose that the condition (\ref{sparse*1}) holds true.
  Then, for all $\a \in (0,\e)$ and $\b \in (0,\8)$, where $\e$ is from (\ref{green*dec}),
  there exists $K_1 \Subset \zd$ such that if 
  $K_1 \sub K \Subset \zd$, then 
\bdmn
& & \| H_{V_K,\lm} \|_{\cB(\ell^2)} \le \b,
\label{H*W*norm}\\
& & H_{V_K,\lm} \in  \cB(\ell^{\8,\a})
\; \; \mbox{with}\; \;
\| H_{V_K,\lm}\|_{\cB(\ell^{\8,\a})} \le \b. 
\label{H*W*inv}
\edmn
\end{lemma}

\noindent Proof: 
Set $\del \deq \e-\a>0$.
  Then, by (\ref{sparse*1}), there exists a $K_1 \Subset \zd$ such that 
  \bdnl{H*W*K}
  \sup_{x \in \zd \bsh K_1}\sum_{y \in \zd  \atop y \neq x}\sqrt{V(x)V(y)}\exp (-\del |x-y|)
  \le {\b \over |\lm|C},
  \edn
  where the constant $C$ is from (\ref{green*dec}). 
  Let $K_1 \sub K \Subset \zd$.
  By (\ref{green*dec}) and (\ref{H*W*ker}), 
  \bdnl{H*W*ker*2}
H_{V_K,\lm}(x,y)
  \le |\lm|C \sqrt{V_K(x)V_K(y)}\exp (-\e |x-y|){\bf 1}_{x \neq y}. 
  \edn
  Since $H_{V_K,\lm}:\ell^2 \ra \ell^2$ is symmetric, we have
  by a standard estimate that 
  \bdnn 
  \| H_{V_K,\lm} \|
  &\le& \sup_{x \in \zd}\sum_{y \in \zd}|H_{V_K,\lm}(x,y)|\\
  &\st{\scriptsize \mbox{(\ref{H*W*ker*2})}}{\le}&
  |\lm|C\sup_{x \in \zd \bsh K}
  \sum_{y \in \zd  \atop y \neq x}\sqrt{V(x)V(y)}\exp (-\e |x-y|)
  \st{\scriptsize \mbox{(\ref{H*W*K})}}{\le} \b,
  \ednn
  which shows (\ref{H*W*norm}).

  Suppose that $u \in \ell^{\8,\a}$ and $x \in \zd$.  
  Then, noting that $|x| \le |x-y|+|y|$, 
  \bdnn 
  \lefteqn{\exp (\a |x|)|(H_{V_K,\lm}u)(x) |} \\
  &\le& \exp (\a |x|)\sum_{y \in \zd}|H_{V_K,\lm}(x,y)||u(y)|\\
  &\st{\scriptsize \mbox{(\ref{H*W*ker*2})}}{\le}&
       {\bf 1}_{\zd \bsh K} (x)|\lm|C
       \sum_{y \in \zd  \atop y \neq x}\sqrt{V(x)V(y)}\exp (-\del |x-y|)\exp (\a |y|)|u(y)|\\
  &\le&
  {\bf 1}_{\zd \bsh K} (x)|\lm|C\|u\|_{\8,\a}
  \sum_{y \in \zd  \atop y \neq x}\sqrt{V(x)V(y)}\exp (-\del |x-y|)
  \st{\scriptsize \mbox{(\ref{H*W*K})}}{\le}
  \b \|u\|_{\8,\a},
  \ednn
  which shows (\ref{H*W*inv}). 
  \hfill $\Box$

  \vvs 
  Combining \Lem{1-gm*W} and \Lem{H*W}, we obtain the following 
  \Lemma{G*W}
  Suppose that (\ref{sparse*1}) holds true.
  Then, for $\lm \in \C \bsh \s_{\rm ess}(P_V)$ and $\a \in (0,\e)$, 
  where $\e$ is from (\ref{green*dec}), there exists $K \Subset \zd$ such that
$1 \in \rh (G_{V_K,\lm})$ and that $(1-G_{V_K,\lm})^{-1} \in \cB(\ell^{\8,\a})$.
\end{lemma}

\noindent Proof: 
We have $\lm \not\in \s_{\rm ess}(P_V)=\s (P) \cup \Lm_V$,
by the choice of $\lm$ and \Thm{ess}.  Thus, we may apply \Lem{1-gm*W}
and take $\e_0 \in (0,1]$ and $K_0 \Subset \zd$ from there.
We then apply \Lem{H*W} with $\b<\e_0$.
As a consequence, we can find $K_i \Subset \zd$ ($i=0,1$) with which
(\ref{1-gm*W}), (\ref{H*W*norm}) and (\ref{H*W*inv}) hold true, where $\b<\e_0$.
We now suppose that $K_0 \cup K_1 \sub K \Subset \zd$. 
Then, it follows from (\ref{1-gm*W}) that the operator $1-\gm V_K$ has 
its inverse in $\cB(\ell^2)$ with the norm at most $1/\e_0$.
Then, setting
$$
R_K=(1-\gm V_K)^{-1}H_{V_K,\lm}
$$
for simplicity, we have 
\bdnl{R*K}
\|R_K\|_{\cB(\ell^2)} \le \|(1-\gm V_K)^{-1}\|_{\cB(\ell^2)} \|H_{V_K,\lm}\|_{\cB(\ell^2)}
\le \e_0^{-1} \cdot \b<1,  
\edn 
and therefore, the following Neumann series converges:
\bdnl{Neumann}
(1-R_K)^{-1}=\sum_{n=0}^\8R_K^n. 
\edn 
Since 
$$
1-G_{V_K,\lm}
=1-\gm V_K-H_{V_K,\lm}
=(1-\gm V_K)(1-R_K),
$$
the operator $1-G_{V_K,\lm}$ is invertible in $\cB(\ell^2)$,
and herefore, $1 \in \rh (G_{V_K,\lm})$. 

By repeating the same argument as above, with
$\cB (\ell^2)$ replaced by $\cB (\ell^{\8,\a})$,
we see that $(1-G_{V_K,\lm})^{-1} \in \cB (\ell^{\8,\a})$.
Indeed, it is clear that $1-\gm V_K$ has the inverse in $\cB (\ell^{\8,\a})$
with the norm at most $\e_0^{-1}$.
Moreover, by (\ref{H*W*inv}), $\|H_{V_K,\lm}\|_{\cB (\ell^{\8,\a})} \le \b <\e_0$.
Consequently, we see that $\| R_K \|_{\cB(\ell^{\8,\a})} < 1$ and hence that
the Neumann series (\ref{Neumann}) converges in $\cB(\ell^{\8,\a})$.
This implies that $(1-G_{V_K,\lm})^{-1} \in \cB(\ell^{\8,\a})$. 
\hfill $\Box$

\vvs
As is discussed in section \ref{decay:outline}, \Thm{decay} follows from the following 
  \Lemma{P*W}
  Suppose that (\ref{sparse*1}) holds true.
  Then, for $\lm \in \C \bsh \s_{\rm ess}(P_V)$ and $\a \in (0,\e)$
  where $\e$ is from (\ref{green*dec}),  
there exists $K \Subset \zd$ such that
$\lm \in \rh (P_{V_K})$ and that $(\lm-P_{V_K})^{-1} \in \cB (\ell^{\8,\a})$.
\end{lemma}

\noindent Proof: 
The first half of the lemma follows from \Lem{H*W}, since, 
$\lm \in \rh (P_{V_K})$ if and only if $1 \in \rh (G_{V_K,\lm})$
by \Lem{BS}. As for the latter half,
by applying (\ref{eq:BS}) with $V$ there replaced by $V_K$,
we have
$$
(\lm - P_{V_K})^{-1}=G_\lm+G_\lm V_K^{1/2}(1-G_{V_K,\lm})^{-1}V_K^{1/2}PG_\lm. 
$$
It is clear that $V_K^{1/2}$ and $P$ belong to $\cB (\ell^{\8,\a})$.
It is also clear from the proof of \Lem{H*W} that $G_\lm \in \cB (\ell^{\8,\a})$.
Moreover, $(1-G_{V_K,\lm})^{-1} \in \cB (\ell^{\8,\a})$ by \Lem{G*W}.
Combining these, we conclude that $(\lm - P_{V_K})^{-1} \in \cB (\ell^{\8,\a})$. 
\hfill $\Box$
\SSC{Proof of \Thm{gap}, \Cor{gap} and \Prop{cond*gap}}
\subsection{Outline}
{\bf Proof of \Thm{gap}}: 
We first prove that
\bdnl{1<dimKer<8}
1 \le \dim \Ker (r(P_V)-P_V)<\8
\; \; \mbox{and}\; \;
r(P_V) \not\in \ov{\s (P_V) \bsh \{r(P_V)\}}.
\edn
For this purpose, we prove via
a ``lower bound by a delta potential'' (\Lem{u0})
that there exists $u_0 \in \ell^2$
such that
\bdnl{u0}
\mbox{$\|u_0\|_V=1$ and
  $\lan P_Vu_0, u_0 \ran_V >\max \s_{\rm ess}(P_V)$},
\edn 
which implies (\ref{1<dimKer<8}) as follows. 
It follows from (\ref{u0}) that 
$$ 
r(P_V) =\sup_{\|u\|_V=1}\lan P_Vu, u \ran_V \ge  \lan P_Vu_0, u_0 \ran_V>\max \s_{\rm ess}(P_V).
$$
Thus, $r(P_V) \not\in \s_{\rm ess}(P_V)$, or equivalently,
$$
\dim (r(P_V)-P_V)<\8\; \; \mbox{and}\; \;
r(P_V) \not\in \ov{\s (P_V) \bsh \{r(P_V)\}},
$$
(cf. (\ref{sig*ess}) and (\ref{LSF4})).
Since $r(P_V) \in \s(P_V)$, $r(P_V)$ is an
eigenvalue and therefore, we obtain (\ref{1<dimKer<8}).

Now that we have established (\ref{1<dimKer<8}),
the rest of the proposition follows immediately
from an extension of the Perron-Frobenius theorem to infinite dimensional setting
(See \Lem{202} below). Indeed,
$$
P_V(x,y)=(1+V(x))P(x,y) \ge P(x,y)
\; \; \mbox{for all $x,y \in \zd$},
$$
which shows that $P_V$ is positive and irreducible. Moreover,
by \Lem{201} below, 
the operator norm of $P_V:\ell^2_V \ra \ell^2_V$ equals to $r(P_V)$.
\hfill $\Box$

\vvs
    \noindent {\bf Proof of \Cor{gap}}
    $P_V$ is positive and irreducible. 
Moreover, $P_V$ is bipartite w.r.t. $J \in \{-1,1\}^{\zd}$ 
if and only if $P$ is bipartite w.r.t. $J$.
We can pass from \Thm{gap} to \Cor{gap} by applying
\Lem{203} below to $T=P_V$.  
    \hfill $\Box$

    \vvs
    \noindent {\bf Proof of \Prop{cond*gap}}: 
              {\it Case1}: Suppose that $p$ vanishes on $A$.
              We prove that $P$ is bipartite w.r.t. $J={\bf 1}_A-{\bf 1}_{A^{\sf c}}$. 
Indeed, it follows from assumptions that ${\supp }[p] \sub A^{\sf c}$,
and therefore by definition of $A$ that
$(x,y) \not\in A^2 \cup (A^{\sf c})^2$ if $y-x \in {\supp }[p]$.
This implies that $P$ is bipartite w.r.t. $J={\bf 1}_A-{\bf 1}_{A^{\sf c}}$. 

{\it Case2}: Suppose that (\ref{K*bsh*0}) holds.
We prove that $-r(P_V)<\ell (P_V)$ for
any nonnegative $V \in \ell^\8(\zd)$. For this purpose, we
use the following inequality, cf. \Lem{edge} below.
\bdnl{edge2}
-r(P_V) +2\ell ({\bf 1}_A P {\bf 1}_A) \le \ell (P_V). 
\edn
On the other hand, we have 
$$
    \ell ({\bf 1}_A P {\bf 1}_A) = 
    \min_{\tht \in [-\pi,\pi]^d}\sum_{x \in A}p(x)\exp (\bfi \tht \cdot x)
     \ge  p(0)-\sum_{x \in A \bsh \{0\}}p(x)>0. 
$$
Then, it follows from  (\ref{edge2}) 
that $-r(P_V) < \ell (P_V)$.
    \hfill $\Box$
\subsection{Lemmas}
We start by explaining the ``lower bound by a delta potential'' referred
to in the proof of \Thm{gap}. 
\Lemma{u0}
\bds
\item[a)]
  Suppose that $v_0>0$. Suppose also that
  $\lm \in (1,\8)$, $v>0$, $g_\lm (0)=1+v^{-1}$, and $V(\zd) \cap [v,\8) \neq \epty$. 
    Then, there exists a $u_\lm \in \ell^2$ such that
  \bdnl{u0a}
  \mbox{$\| u_\lm \|_V=1$, and $\lan P_Vu_\lm, u_\lm \ran_V >\lm$}.
  \edn 
  \item[b)]
    Under hypothesis of \Thm{gap},
    there exists a $u_0 \in \ell^2$ such that (\ref{u0}) holds.
    \eds
\end{lemma}
Proof: 
a) We take $c \in \zd$ such that $V(c) \ge v$. 
Then, by (\ref{irred}), $\vp_\lm (x)\deq G_\lm (x,c)>0$ for all $x \in \zd$.
We will show that the function $u_\lm \deq \vp_\lm/\| \vp_\lm \|_V$ satisfies (\ref{u0a}).
To see this, we make an auxiliary use of the potential $U=v \del_c$
to give a lower bound for $V$. 
We first observe that
\bdnl{vp*ev}
P_Uu_\lm =\lm u_\lm. 
\edn
Indeed, noting the equality $PG_\lm=\lm G_\lm-1$, we have for all $x \in \zd$ that
\bdnn
P_U\vp_\lm (x)
&=& (1+v \del_c(x))PG_\lm \del_c (x)\\
&=& (1+v \del_c(x))(\lm G_\lm (x,c)-\del_c (x))\\
&=& \lm G_\lm (x,c)-\del_c (x)+v(\lm G_\lm (c,c)-1)\del_c (x)\\
&=&\lm G_\lm (x,c)=\lm \vp_\lm (x),
\ednn
which implies (\ref{vp*ev}).

Now, we use (\ref{vp*ev}) to see (\ref{u0a}) as follows. 
It follows from $v \le V(c)$ that $U(x) \le V(x)$ for all $x \in \zd$.
Moreover, since $v_0>0$, $V(x)>0$ for infinitely many $x$'s, and hence 
$U(x) < V(x)$
for at least an $x$ (in fact, for infinitely many $x$'s) . 
These, together with the strict positivity of $u_\lm$ implies that 
\bdnl{||u||<1}
1=\| u_\lm \|_V=\|(1+V)^{-1/2}u_\lm\|<\|(1+U)^{-1/2}u_\lm\|=\| u_\lm \|_U.
\edn
Moreover, since $\lan Pu_\lm,u_\lm \ran >0$,
$$
\lan P_Vu_\lm,u_\lm \ran_V=\lan Pu_\lm,u_\lm \ran
\st{\scriptsize \mbox{(\ref{||u||<1})}}{>}
   {\lan Pu_\lm,u_\lm \ran \over \| u_\lm \|_U}
={\lan P_Uu_\lm,u_\lm \ran_U \over \| u_\lm \|_U}
\st{\scriptsize \mbox{(\ref{vp*ev})}}{=}
\lm.
$$
b) Here, by \Cor{ess},
$\lm_0=\max \s_{\rm ess}(P_V)$ and $v_0$ are related as
 $$
 g_{\lm_0}(0)=1+v_0^{-1}.
 $$
 Therefore, part b) of the present lemma follows from part a). 
 \hfill $\Box$

 \vvs
 We present a series of  lemmas in the following more general framework.  
 Let $(S,\cA,m)$ be a $\s$-finite measure space.
 We denote  the norm and the the inner product of the Hilbert space
 $L^2(m)$ resectively by
 $\| \cdot \|$ and $\lan\cdot,\cdot\ran$. The totality of
 $f \in L^2(m)$ such that $f \ge 0$, $m$-a.e. is denoted by $L^2_+(m)$. 
An operator $T \in \cB (L^2(m))$ 
 is said to be {\it positive} if $Tf \in L^2_+(m)$ for all $f \in L^2_+(m)$.
\Lemma{201}
Suppose that $T \in \cB (L^2(m))$ is positive. 
Then, for all $f \in L^2(m)$, 
\bdnl{201*1}
T|f|-|Tf| \in L^2_+(m),
\; \; \mbox{and thus,}\; \;
|\lan f,Tf \ran | \le \lan |f|,T|f| \ran.
\edn
In particular, if $T$ is positive and self-adjoint, then, 
$$
|\ell (T)| \le r (T)=\|T\|.
$$
\end{lemma}
  Proof: 
  Let $f \in L^2(m)$ be arbitrary.
  Since $T$ maps a real-valued function to a real-valued function, we have 
  \bdnl{TRe=ReT}
  T (\Re f)=\Re (Tf).
  \edn
  Moreover, for each $z \in \C$, 
  \bdnl{|z|=}
|z|=\sup_{\tht \in \R}\Re (\exp (\bfi \tht) z)=\sup_{\tht \in \Q}\Re (\exp (\bfi \tht) z). 
\edn
Let $x \in S$ and $\tht \in \Q$ be fixed.
Then, by applying (\ref{|z|=}) to $z=f(x)$, 
and to $z=Tf(x)$,
we see that 
\bdmn
g_\tht(x) & \deq & |f(x)| -\Re (\exp (\bfi \tht) f(x))
\st{\scriptsize \mbox{(\ref{|z|=})}}{\ge} 0,
\label{|z|=*1}\\
  T|f|(x)-|Tf(x)|
  &\st{\scriptsize \mbox{(\ref{|z|=})}}{=}&
  T|f|(x)-\sup_{\tht \in \Q}\Re (\exp (\bfi \tht) Tf(x))\nn \\
  &=&
  \inf_{\tht \in \Q}\lef( T|f|(x)-\Re (\exp (\bfi \tht) Tf(x))\ri)\nn \\ 
  &\st{\scriptsize \mbox{(\ref{TRe=ReT})}}{=}&
\inf_{\tht \in \Q}Tg_\tht(x).
\label{|z|=*2}
  \edmn
  By (\ref{|z|=*1}), we have $Tg_\tht \in L^2_+(m)$.
  Then, combining this with (\ref{|z|=*2}), we have $T|f|-|Tf| \in L^2_+(m)$. 
  In particular, if $T$ is positive and self-adjoint, then,
  $$
  -\lan f,T f \ran
  \st{\scriptsize \mbox{(\ref{201*1})}}{\le} \lan |f|, T|f| \ran \le r(T)\|f\|^2. 
  $$
Taking the supremum of the left-hand side of the above inequality
over all normalized functions $f$ in $L^2 (m)$,
we obtain $-\ell (T)  \le r (T)$, which, together with
the obvious inequality $\ell (T)  \le r (T)$, implies that
$|\ell (T)| \le r (T)$. As a consequence, we obtain the
equality $r (T)=\|T\|$.
\hfill $\Box$

\vvs 
Let $T \in \cB (L^2(m))$. $T$ is said to be {\it irreducible}
if, for all functions $f,g \in L^2_+(m)$ which are not
identically zero $m$-a.e., there exists $n \in \N$ such that
$\lan f,T^ng \ran > 0$.
If there exists a measurable function $J:S \ra \{-1,1\}$
such that $JT=-TJ$, i.e., $JTf=-T(Jf)$ for all $f \in L^2(m)$,
then, $T$ is said to be {\it bipartite} and $J$ is called the
{\it sign} of $T$. 
$T$ is simply said to be bipartite if $T$ is bipartite w.r.t. some
sign $J$. 

We consider the following conditions for
an operator $T \in \cB (L^2(m))$, in connection with
its irreducibility and bipartiteness. 
\bdmn
& & \Ker (\| T \|+T) \neq \{0\}.
\label{psi*202*0}\\
& & \Ker (\| T \|-T) \neq \{0\}.
\label{vp*202*0}\\
& & \mbox{There exists a normalized function $\vp \in L^2(m)$ such that} \nn\\
& & \vp > 0\; \mbox{a.e. and}\;\Ker (\| T \|-T)=\C \vp.
\label{vp*202*1}
\edmn
\Lemma{202}
Suppose that $T \in \cB (L^2(m))$ is positive. Then, 
\bdnl{202*a} 
\mbox{(\ref{vp*202*1}) $\Llra$ (\ref{vp*202*0}) and $T$ is irreducible}.
\edn
Suppose in addition that $T$ is irreducible. Then,
\bdnl{202*b}
\mbox{(\ref{psi*202*0}) $\Llra$ (\ref{vp*202*1}) and that $T$ is bipartite.}
\edn
Moreover, the converse part of the above equivalence entails the following
relation. If $T$ is bipartite w.r.t. the sign $J$,
then, (\ref{psi*202*0}) holds with  
\bdnl{psi*202*1}
  \Ker (\| T \|+T)=\C J\vp.
  \edn 
   \end{lemma}
Proof: Let us make a preliminary observation which applies 
to all $T \in \cB(L^2(m))$.  
Suppose that $u \in L^2(m)$ is normalized and that $\lan u,Tu \ran =\s \|T\|$,
where $\s$ is either $1$ or $-1$. Then,
\bdnl{202*obs}
Tu =T^*u=\s \|T\|u.
\edn
Indeed, since $\lan u,Tu \ran =\s \|T\| \in \R$, 
 $$
  \|\s \|T\|u -Tu\|^2
  =
  \|T\|^2-2\s \|T\|\lan u, Tu \ran  +\| Tu \|^2
  =-\|T\|^2+\| Tu \|^2 \le 0.
  $$
  Thus, $Tu =\s \|T\|u$.
  On the other hand, $\lan u,Tu \ran =\s \|T\|$ implies that $\lan u,T^*u \ran =\s \|T^*\|$.  
Hence by letting $T^*$ play the role of $T$ above, we obtain $T^*u=\s \|T\|u$ as well. 
\\ (\ref{202*a}): This equivalence is due to \cite[p.202, Theorem XIII.43]{RS78}.
\\  (\ref{202*b}) ($\Ra$):
Suppose that $\psi \in L^2(m)$ is a nonzero function which satisfies the equation 
$T\psi=-\|T\|\psi$. 
Since the same equation is satisfied by the real and imaginary parts of $\psi$,
we may assume that $\psi$ is real-valued.
  Moreover, by normalization, we may assume that
  $\|\psi\|=1$. Therefore, we see from the preliminary
  observation (\ref{202*obs}) that
  \bdnl{202*obs*1}
T\psi =T^*\psi =-\|T\|\psi. 
  \edn 
  On the other hand, applying (\ref{201*1}) to $f=\psi$,
  we obtain
  $$
  \|T\| =-\lan \psi, T\psi \ran \le \lan |\psi|, T|\psi| \ran 
  \le \|T\|.
  $$ 
  and hence, $\lan |\psi|, T|\psi| \ran =\|T\|$. Therefore, we see from the preliminary
  observation (\ref{202*obs}) that $ T|\psi|=\|T\||\psi|$, 
  Thus, we have proved (\ref{vp*202*0}), which, together with the irreducibility of $T$,
  implies (\ref{vp*202*1}).
  Then, it follows from (\ref{vp*202*1}) that $|\psi|=\vp$, and by
  applying (\ref{202*obs}), we have
  \bdnl{202*obs*2}
T\vp =T^*\vp =\|T\|\vp. 
  \edn 
  Since $|\psi|=\vp$, the function $J \deq \psi /\vp$ is defined a.e. and takes
  values in $\{-1,1\}$ there. We will prove that $T=-JTJ$. By linearlity,
  it is enough to prove that $g \deq Tf+JTJf=0$ a.e. for all $f \in L^2_+(m)$.
  Since $f \pm Jf \in L^2_+(m)$, we have $Tf \pm TJf \in L^2_+(m)$,
  and hence $g \in L^2_+(m)$. Therefore, it is enough to show
  that $\lan \vp, g \ran =0$,
  which can be done by noting that $\vp J=\psi$, and that $\psi J=\vp$ as follows.
  \bdnn
  \lan \vp, g \ran
  &=&\lan \vp, Tf \ran+\lan \psi, TJf \ran
  =\lan T^*\vp, f \ran+\lan T^*\psi, Jf \ran\\
  &\st{\scriptsize \mbox{(\ref{202*obs*1}). (\ref{202*obs*2})}}{=}
  &\|T\| \lan \vp, f \ran-\|T\|\lan \psi, Jf \ran
= \|T\| \lan \vp, f \ran-\|T\|\lan \vp, f \ran=0. 
\ednn
\\  (\ref{202*b}) ($\La$):
By the relation $TJ=-JT$, the eigenspaces $\Ker (\|T\| \pm T)$ are mapped to each other bijectively by $J$. 
  Thus, (\ref{vp*202*1}) implies (\ref{psi*202*1}). 
\hfill $\Box$
  
\vvs 
\Lemma{203}
Suppose that an operator $T \in \cB (L^2(m))$ is positive, self-adjoint, 
and satisfies (\ref{vp*202*1}) and
\bdnl{SGT}
r (T)  >  \sup\{\lm \; ; \; \lm \in \s (T),\; \lm \neq r (T)\}. 
\edn 
Then,
\bds
\item[a)] The following conditions are equivalent.
\bdmn
& &\mbox{$r (T)  >  \sup\{|\lm| \; ; \; \lm \in \s (T),\; |\lm| \neq r (T)\}$.}
\label{203*1}\\
& &\mbox{Either $-r (T) <\ell (T)$ or $T$ is bipartite.}
\label{203*2}
\edmn
\item[b)] Suppose that one of the conditions (\ref{203*1}) and (\ref{203*2}) 
  holds true (therefore that both do). 
Then, there exists $\e \in (0,1)$ such that
\bdnl{proj2*204}
\| r(T)^{-n}T^n (f -\Pi f)\| \le \e^n\|f -\Pi f\|
\; \; \;
\mbox{for all $f \in L^2(m)$ and $n \in \N$},
\edn
where, with $\vp$ from (\ref{vp*202*1}), 
\bdnl{proj*204}
\Pi f \deq 
\lef\{\barray{ll}
\lan f, \vp \ran \vp, 
& \mbox{if $-r (T)<\ell (T)$},\\
\lan f, \vp \ran \vp+\lan f, J\vp \ran J\vp,
& \mbox{if $T$ is bipartite w.r.t. the sign $J$}.
\earray \ri. 
\edn
  \eds
\end{lemma}
  Proof: 
  a) (\ref{203*1}) $\Ra$ (\ref{203*2}): 
  Suppose that (\ref{203*2}) fails. Then, it follows from
  \Lem{202} a) that $\ell (T)=-r(T)$ is not an
  eigenvalue, and therefore by (\ref{LSF4}), $-r(T) \in \ov{\s(T) \bsh\{ -r(T)\}}$.
  Hence (\ref{203*1}) fails. 
  \\ (\ref{203*2}) $\Ra$ (\ref{203*1}): 
  If $-r (T)<\ell (T)$, then (\ref{203*1}) holds, since
  \bdnl{sig*sub*203}
\s (T)  \cap   [-r (T), \ell (T))=\epty.
       \edn
  On the other hand, if $T$ is
  bipartite, then, (\ref{203*1}) follows from (\ref{SGT}) and the 
  symmety of $\s(T)$ with respect to the origin. 
  \\ b) Let $T =\int_{\s (T)}\lm d E_\lm$ denote the spectral decomposition of
the self-adjoint operator $T$ on $\ell^2(m)$. 
    We first verify that
    \bdnl{Pi=*204}
\Pi = \int_{\s (T) \atop |\lm|=r(T)}dE_\lm. 
\edn 
We treat each of the cases in condition (\ref{203*2}) separately. 
Let us temporarily denote the orthogonal projection on
the right-hand side of (\ref{Pi=*204}) by $\Pi^{\rm RHS}$. 

 {\it Case1}: Suppose that $-r (T) <\ell (T)$. 
 Then, it follows from (\ref{sig*sub*203}) that
 $$
 \{\lm \in \s (T)\; ; \; |\lm|=r(T)\}=\{r(T)\}. 
 $$
We see from this and (\ref{vp*202*1}) that 
$$
\Ran (\Pi^{\rm RHS})=\Ker (r(T)-T)=\C \vp,
$$
which implies (\ref{Pi=*204}). 

{\it Case2}: Suppoes that $T$ is bipartite the sign $J$. 
Then, it follows from the symmetry of $\s(T)$ with respect to the origin
that
$$
 \{\lm \in \s (T)\; ; \; |\lm|=r (T)\}=\{\pm r (T)\}. 
 $$
 Putting this, (\ref{vp*202*1}) and (\ref{psi*202*1}) together, we have 
$$
\Ran (\Pi^{\rm RHS})
= \Ker (r (T)-T) \oplus \Ker (-r (T)-T)
=\C \vp \oplus \C J\vp,
$$
which implies (\ref{Pi=*204}).

Now, (\ref{proj2*204}) follows easily from condition (\ref{203*1}).  
Let
$$
r=\sup\{|\lm| \; ; \; \lm \in \s (T),\; |\lm| \neq r (T)\}
$$
and let $f^\perp=f-\Pi f$.
Then, it follows from (\ref{Pi=*204}) that 
$$
\| T^nf^\perp\|_V^2
=\int_{\s (T) \atop |\lm|\le r}\lm^{2n}d\lan E_\lm f^\perp, f^\perp \ran
\le r^{2n}\|f^\perp\|^2,
$$
which proves (\ref{proj2*204}). 
\hfill $\Box$
\Lemma{edge}
Suppose that $p:\zd \ra [0,\8)$ is a transition probability of a
  symmetric random walk,
  and that $A=\{x \in \zd\; ; \; \sum_{\a \in I}x_\a \in 2\Z\}$
  with $\epty \neq I \sub \{1,\ldots,d\}$. 
  Then,   
  \bdnl{left*edge}
  -r(P_V) +2\ell ({\bf 1}_A P {\bf 1}_A) \le \ell (P_V).
\edn
  \end{lemma}
  Proof: 
We already know from \Lem{201} that $-r(P_V) \le \ell (P_V)$. 
Thus, we may assume that $\ell ({\bf 1}_A P {\bf 1}_A) \ge 0$. Since
  $$
-\ell (P_V) =\sup \{ -\lan P_Vu,u\ran_V \; ; \; \|u\|_V=1\},
$$
it is enough to show that 
\bdnl{edge3}
-\lan P_Vu,u\ran_V \le (-2\ell ({\bf 1}_AP{\bf 1}_A)+r(P_V))\|u\|^2_V
\; \; \mbox{for all $u \in \ell^2$}.
\edn
To see this, let $u \in \ell^2$ be arbitrary and $J={\bf 1}_A-{\bf 1}_{A^{\sf c}}$.
We first verify that
  \bdnl{edge4}
-\lan Pu,u\ran \le -2\ell ({\bf 1}_A P {\bf 1}_A)\|u\|_2^2+\lan PJu,Ju\ran.
\edn
Let
$P_+={\bf 1}_A P {\bf 1}_A +{\bf 1}_{A^{\sf c}}P{\bf 1}_{A^{\sf c}}$ 
and
$P_{-}={\bf 1}_A P{\bf 1}_{A^{\sf c}}+{\bf 1}_{A^{\sf c}}P{\bf 1}_A$. 
Then, 
\bdnl{edge5}
P_{\pm}J=\pm JP_{\pm},
\edn
as is easily be seen. 
Moreover, for fixed $c \in A^{\sf c}$,
the map $x \mapsto x+c$ ($A$ $\lra$ $A^{\sf c}$) is a bijection.
Hence
\bdmn
\lan P_+u,u\ran
&=& \lef(\sum_{x,y \in A}+\sum_{x,y \not\in A}\ri)p(y-x)u(x)u(y)^*\nn \\
&=& \sum_{x,y \in A}p(y-x)u(x)u(y)^* \nn \\
& & +\sum_{x',y' \in A}p(y'-x')u(x'+c)u(y'+c)^* \nn \\
& \ge & \ell ({\bf 1}_A P {\bf 1}_A)\sum_{x \in A} \lef( |u(x)|^2+|u(x+c)|^2\ri)
=\ell ({\bf 1}_A P {\bf 1}_A)\|u\|^2. 
\label{edge6}
\edmn
  Thus, we obtain (\ref{edge4}) as follows.
  \bdnn
  \lan PJu,Ju\ran
  &=&   \lan P_+Ju,Ju\ran+\lan P_-Ju,Ju\ran\\
  &\st{\scriptsize \mbox{(\ref{edge5})}}{=}&
  \lan JP_+u,Ju\ran-\lan JP_-u,Ju\ran\\
  &=&   \lan P_+u,u\ran-\lan P_-u,u\ran
  =   2\lan P_+u,u\ran-\lan Pu,u\ran \\
  & \st{\scriptsize \mbox{(\ref{edge6})}}{\ge} &
  2 \ell ({\bf 1}_A P {\bf 1}_A)\|u\|_2^2-\lan Pu,u\ran.
  \ednn 
    We then use (\ref{edge4}) to prove (\ref{edge3}) as follows.
    Note that $\ell ({\bf 1}_A P {\bf 1}_A)\|u\|^2\ge \ell ({\bf 1}_A P {\bf 1}_A)\|u\|_V^2 $,
    since $\|u\| \ge \|u\|_V$ and we have assumed that $\ell ({\bf 1}_A P {\bf 1}_A) \ge 0$.
    Therefore, 
    \bdnn
 -\lan P_Vu,u\ran_V 
 &=&-\lan Pu,u\ran\\ 
 & \st{\scriptsize \mbox{(\ref{edge4})}}{\le} &
 -2\ell ({\bf 1}_A P {\bf 1}_A)\|u\|^2+\lan PJu,Ju\ran\\
  & = &
  -2\ell ({\bf 1}_A P {\bf 1}_A)\|u\|^2+\lan P_VJu,Ju\ran_V\\
  &\le&
  (-2\ell ({\bf 1}_A P {\bf 1}_A)+r(P_V))\|u\|^2_V. 
  \ednn
\hfill $\Box$

\small
    \noindent {\bf Acknowledgments}
The second author was supported by JSPS KAKENHI
Grant Number 25400136


\end{document}